\documentclass[preprint,12pt]{elsarticle}

\usepackage{pgf}
\usepackage{graphicx,color}
\usepackage{bm}
\usepackage{multirow}

\newcommand{\ez}{E}
\newcommand{\pz}{P}

\newcommand{\eqd}{\stackrel{D}{=}}

\usepackage{amssymb,amsfonts,amsmath}
\usepackage{amsthm}

\newtheorem{thm}{Theorem}
\newtheorem{lem}[thm]{Lemma}

\journal{Mathematical Biosciences}

\begin{document}

\begin{frontmatter}



\title{Stochastic epidemic models featuring contact tracing with delays}


\author{Frank G.~Ball}
\ead{frank.ball@nottingham.ac.uk}
\author{Edward S.~Knock\corref{cor1}}
\ead{edwardknock@gmail.com}
\author{Philip D.~O'Neill}
\ead{philip.oneill@nottingham.ac.uk}
\cortext[cor1]{Corresponding author.}

\address{School of Mathematical Sciences, University of Nottingham, University Park, Nottingham NG7 2RD, UK.}

\begin{abstract}
This paper is concerned with a stochastic model for the spread of an SEIR (susceptible $\rightarrow$ exposed (=latent) $\rightarrow$ infective $\rightarrow$ removed) epidemic with a contact tracing scheme, in which removed individuals may name some of their infectious contacts, who are then removed if they have not been already after some tracing delay. The epidemic is analysed via an approximating, modified birth-death process, for which a type-reproduction number is derived in terms of unnamed individuals, that is shown to be infinite when the contact rate is sufficiently large. We obtain explicit results under the assumption of either constant or exponentially distributed infectious periods, including the epidemic extinction probability in the former case. Numerical illustrations show that, while the distributions of latent periods and delays have an effect on the spread of the epidemic, the assumption of whether the delays experienced by individuals infected by the same individual are of the same or independent length makes little difference.
\end{abstract}


\end{frontmatter}



\section{Introduction}
\label{intro}
This paper is concerned with a stochastic model for the spread of an infectious disease amongst a homogeneously mixing population. The model features a contact tracing scheme in which individuals can name their past infectious contacts upon diagnosis, with traced contacts being prevented from further infecting others. Such a model was examined in Ball {\em et al.}~\cite{ball-knock-oneill}, wherein it was assumed that at the end of an individual's infectious period they name, at random, people they had infected, with named undiagnosed individuals being traced and prevented from making further potentially infectious contacts (in practice, via isolation or vaccination, for instance) and with some probability a traced individual then also names, at random, their infectious contacts. It was assumed that the tracing process was instantaneous, a simplifying assumption that is unrealistic, particularly when there can exist a long tracing chain in which an individual is diagnosed, names a contact, who is then traced and in turn names a contact, who is then traced and names a contact and so on; such a chain was considered instantaneous, no matter how long. Further, it was assumed that when an individual is infected they immediately are able to infect others themselves, an assumption which can be unrealistic for diseases which have a non-negligible latent period. As well as increased realism, incorporating a latent period may make it possible for infected individuals to be traced while they are latent and so before they have been able to make any infectious contacts, and as such this can enhance the effect of contact tracing.

In this paper we derive a threshold parameter for a new model which extends that of Ball {\em et al.}~\cite{ball-knock-oneill} by (i) incorporating latent periods and tracing delays; and (ii) assuming that, as was already the case for traced individuals, with some probability an untraced individual is not asked to name their infectious contacts (perhaps because they are asymptomatic). Delays have been incorporated in contact tracing models such as those of Klinkenberg {\em et al.}~\cite{klink-fraser-heest}, who assumed that delay is of fixed length, and Shaban {\em et al.}~\cite{shabanetal}, who assume that, after a delay time beginning at infection, an individual is detected and their contacts traced. As such, these models assume that individuals named by a shared infector experience the same delay. We assume instead that such individuals experience independent delays (as Ross \& Black \cite{ross-black} do), but for comparison also consider the effect delays of the same length would have instead.

Previous stochastic models for contact tracing, such as Klinkenberg {\em et al.}~\cite{klink-fraser-heest} and M\"{u}ller {\em et al.}~\cite{muller-kretzschmar-dietz}, assume that all individuals will eventually be asked to name contacts, while our model assumes there is a probability for each infected individual that they may never be asked to name their contacts. Other tracing models may exploit population structure to target intervention to susceptibles, such as Ross \& Black \cite{ross-black} (wherein antivirals are given to traced infected individuals and their housemates, reducing susceptibility and infectivity) and Shaban {\em et al.}~\cite{shabanetal} (wherein traced neighbours of a diagnosed individual in a network-structured population are vaccinated, rendering susceptibles immune). The tracing in this paper differs from these in that intervention is only targeted towards infected individuals, and further, unlike Shaban {\em et al.}, traced individuals are necessarily infected. Additionally, the models of Shaban {\em et al.}~and Ross \& Black assume that some individuals, once traced, may yet make further contacts, while our model assumes all traced individuals are prevented from making further contacts.

The remainder of the paper is structured as follows. In Section \ref{sec:model-definition} the epidemic with contact tracing model is outlined. A modified birth-death process which approximates the epidemic is described in Section \ref{sec:mbdp}. In Section \ref{sec:const}, a threshold parameter and extinction probability is derived via a two-type branching process for constant infectious periods in the special case where only traced individuals may name contacts. An embedded Galton-Watson process of unnamed individuals for use in more general cases is described in Section \ref{sec:EGWP}, and is then used to determine a threshold parameter for exponential infectious periods in Section \ref{sec:exp}. Numerical illustrations of the theory are presented in Section \ref{sec:numerics} and Section \ref{sec:discussion} contains some further discussion.

\section{Model definition}
\label{sec:model-definition}

We consider an SEIR (Susceptible $\rightarrow$ Exposed $\rightarrow$ Infective $\rightarrow$ Removed) epidemic spreading amongst a homogeneously mixing, closed population of size $N+m$, with a contract tracing scheme applied to reduce spread. At any time, each individual in the population is in one of four states: susceptible, exposed (i.e.~latent), infective or removed. Initially a small number, $m$, of individuals are infectives and the remaining $N$ are susceptible. A susceptible individual becomes a latent individual if he/she makes contact with an infective in a manner described below. A latent individual remains latent for a period of time distributed according to a random variable $T_L$, having an arbitrary but specified distribution (i.e.~no assumption is made about the form of its distribution, but the distribution has to be known), at the end of which he/she becomes infective. An infective individual remains infectious for a period of time distributed according to a random variable $T_I$, having an arbitrary but specified distribution, and then becomes removed. Contacts between two given individuals in the population occur at times given by the points of a homogeneous Poisson process with rate $\lambda/N$. Once removed, an individual no longer plays a part in the epidemic process. The epidemic ends when there are no more latent or infective individuals left in the population. All of the Poisson processes, and the random variables describing latent and infectious periods, are assumed to be mutually independent.

This epidemic incorporates a contact tracing scheme as follows. Upon removal, individuals may be interviewed, as described below. An interviewed individual names each of the individuals they infected, if any, (i.e.~their infectees) independently with probability $p$, and named infectees are removed after a delay period of time distributed according to a random variable $T_D$, having an arbitrary but specified distribution. Generally, we assume that the random variables describing the delay periods of all individuals with the same infector (i.e.~siblings) are mutually independent. However, in some cases it is shown that results hold also when siblings are assumed to have delay periods of the same length. An individual whose removal is a result of contact tracing (and not the natural end of their infectious period) is called traced, otherwise they are untraced. Note that an untraced individual can be either unnamed, or named but their infectious period ends during their infector's infectious period or the associated tracing delay. Untraced and traced individuals are interviewed with probability $\pi_R$ and $\pi_T$, respectively, otherwise all of their infectees are unnamed. The naming process and random variables describing delay periods are assumed independent of the Poisson processes, and random variables describing latent and infectious periods. While the interview probabilities $\pi_R$ and $\pi_T$ could just be considered facets of the tracing mechanism, they could instead be considered the probabilities of an individual being symptomatic when their removal is brought about by the end of their infectious period ($\pi_R$) or their being traced ($\pi_T$), with all symptomatic individuals and no asymptomatic individual being interviewed. Note that under this interpretation, infectives are symptomatic independently, and symptomatic and asymptomatic individuals have the same latent and infectious period distributions and the same infection rate. The most important parameters appearing throughout the paper have been listed in Table \ref{table:parameters}. Note that $R_0$ is defined in Section \ref{sec:const} and $R_U$ and $p_E$ are defined in Section \ref{sec:EGWP}.

\begin{table}[h]
\begin{tabular}{ll}
\hline
Parameter & Description\\
\hline
$N$ & Initial number of susceptibles\\
$m$ & Initial number of infectives\\
$\lambda/N$ & Individual-to-individual contact rate\\
$p$ & Probability that an interviewed individual names a given\\
    & infectee\\
$\pi_R$ & Probability that an untraced individual is interviewed\\
$\pi_T$ & Probability that a traced individual is interviewed\\
$\iota$ & Length of constant infectious period\\
$\gamma$ & Rate parameter for exponentially-distributed infectious period\\
         & (i.e.~$\mbox{mean}=1/\gamma$)\\
$R_0$ & Expected number of offspring of a typical individual in the\\
      & modified birth-death process\\
$R_U$ & Expected number of offspring of a typical individual in the\\
      & embedded Galton-Watson process of unnamed individuals\\
$p_E$ & Extinction probability of the modified birth-death process\\
\hline
\end{tabular}
\caption{List of important parameters}
\label{table:parameters}
\end{table}

\section{Approximating modified birth-death process (MBDP)}
\label{sec:mbdp}

The following modified birth-death process was originally presented in Ball {\em  et al.}~\cite{ball-knock-oneill}, but is extended here to incorporate delays and latent periods.

If the initial number of susceptibles $N$ is large, then during the early stages of the epidemic, there is only a small probability that an infective makes contact with an already-infected individual. Thus we can approximate the early stages of the epidemic by a process in which all of an infective's contacts are made with susceptible individuals. In this approximation, the process of infected individuals follows a modified birth-death process (MBDP). The approximation can be made fully rigorous by considering a sequence of epidemics, indexed by $N$, and modifying the coupling argument of Ball \cite{ball83} (see also Ball \& Donnelly \cite{ball-donnelly}) to prove almost-sure convergence of the epidemic model with tracing to the MBDP as $N\to\infty$.

This MBDP, with births corresponding to new infections and deaths corresponding to removals, is as follows. Individuals give birth at rate $\lambda$ over their active lifetime which has a natural length distributed as $T_I$. An individual's active lifetime begins after a latent period distributed as $T_L$, and ends when they die. When an individual dies they may be interviewed. An interviewed individual names each of its offspring independently and with probability $p$. Named offspring are traced after some delay distributed as $T_D$ (the delays of siblings are independent of each other). An individual who is still alive at the time they are traced then dies immediately. Individuals whose deaths are as a result of being traced are said to die unnaturally, all other individuals die naturally. Individuals who die naturally and unnaturally are interviewed with probability $\pi_R$ and $\pi_T$, respectively. Note that the active lifetime length of an individual who is traced during their latent period is zero and such individuals have no offspring. In the epidemic context, active lifetimes correspond to infectious periods, natural deaths correspond to untraced removals and unnatural deaths correspond to traced removals.

\section{Constant infectious period}
\label{sec:const}

In this section we assume the infectious period has a fixed value $\iota$, i.e.~$T_I\equiv\iota$, and that the distributions of $T_L$ and $T_D$ are arbitrary but specified. In this case it is difficult to make analytical progress when traced individuals can name their offspring as the active lifetimes of such traced individuals are no longer constant. Therefore, we restrict attention to the case where only untraced individuals may be interviewed by setting $\pi_T=0$.

Consider an interviewed individual in the MBDP who has active lifetime $\tau$ and $k$ named offspring. Because the births of these $k$ named offspring occur at the points of a homogeneous Poisson process over the interviewed individual's active lifetime $\tau$, the intervals between birth and being named for these offspring can be described by $k$ independent $\mbox{U}(0,\tau)$ random variables. Hence named offspring co-depend upon their parent's active lifetime. Under the assumptions that lifetimes have fixed natural length $\iota$ and that only untraced individuals may be interviewed, the active lifetime of interviewed individuals must necessarily be of fixed length $\iota$, and so named individuals with the same parent are independent. Hence the behaviour of the MBDP in this case can be analysed via a two-type Galton-Watson process (TTGWP), where the types are unnamed ($U$) and named ($N$)  individuals.

Let $Z_{UU}$ ($Z_{NU}$) and $Z_{UN}$ ($Z_{NN}$) be the number of unnamed offspring and named offspring, respectively, of a typical unnamed (named) individual, and let $m_{UU}=\mbox{E}\left[Z_{UU}\right]$, $m_{UN}=\mbox{E}\left[Z_{UN}\right]$, $m_{NU}=\mbox{E}\left[Z_{NU}\right]$ and $m_{NN}=\mbox{E}\left[Z_{NN}\right]$. For $0\leq s_U,s_N\leq1$, let $f_U\left(s_U,s_N\right)=\mbox{E}\left[s_U^{Z_{UU}}s_N^{Z_{UN}}\right]$ and $f_N\left(s_U,s_N\right)=\mbox{E}\left[s_U^{Z_{NU}}s_N^{Z_{NN}}\right]$. Standard results from branching process theory (see Section 5.5 of Haccou et al.~\cite{haccou}) tell us the TTGWP will die out with probability $1$ if $R_0\leq1$, where $R_0$ is the maximal eigenvalue of the matrix of mean offspring
\begin{align*}
M & = \left[\begin{array}{cc} m_{UU} & m_{UN}\\
			m_{NU} & m_{NN}\end{array}\right],
\end{align*}
and so is given by
\begin{align}
R_0 & = \frac{1}{2}\left(m_{UU}+m_{NN}+\sqrt{\left(m_{UU}-m_{NN}\right)^2+4m_{UN}m_{NU}}\right).
\label{eqn:R0}
\end{align}
If $R_0>1$, then, letting $q_U$ and $q_N$ be the extinction probability given there is one initial individual who is unnamed and named, respectively, $(q_U,q_N)$ is the unique solution in $(0,1)^2$ of
\begin{align*}
q_U & = f_U\left(q_U,q_N\right),\\
q_N & = f_N\left(q_U,q_N\right).
\end{align*}
Naturally we assume that $m$ initial individuals are all unnamed, so the extinction probability is $q_U^m$.

In the remainder of this section we calculate the matrix of mean offspring $M$ and probability generating functions $f_U$ and $f_N$. The mean matrix $M$ is used in the numerical illustrations in Section \ref{sec:numerics} (see Figure \ref{fig:lambdacrit-delay}), and $f_U$ and $f_N$ enable extinction probabilities to be calculated without resorting to simulation when $\pi_T=0$, cf.~Figure \ref{fig:estimatedpE} and Table \ref{table:mutvsinddelays}.

\subsection{Calculation of the mean offspring matrix}
\label{subsec:mean offspring matrix}

The number of offspring of an unnamed individual has a Poisson distribution with mean $\lambda\iota$. With probability $\pi_R$, proportions $p$ and $(1-p)$ of these are expected to be named and unnamed, respectively, otherwise all of the offspring are unnamed. Hence
\begin{align*}
m_{UU} & = \lambda\iota(\pi_R(1-p)+(1-\pi_R)) = \lambda(1-\pi_Rp)\iota,\\
m_{UN} & = \lambda\pi_Rp\iota.
\end{align*}

\begin{figure}[h!]
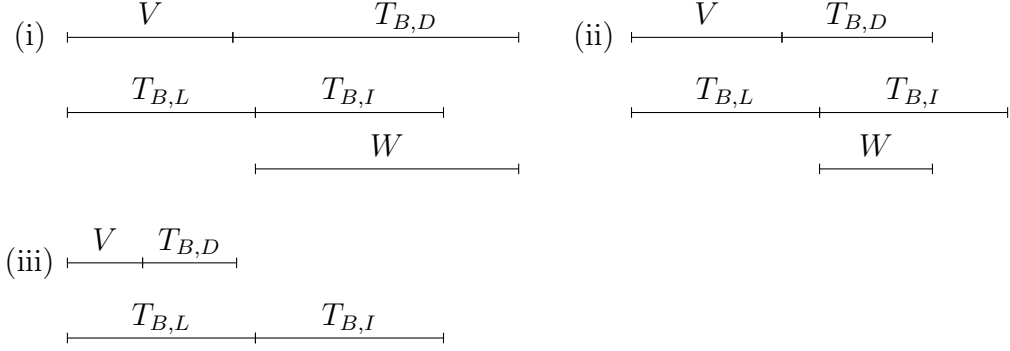

\begin{pgfpicture}{-1cm}{8cm}{13.5cm}{12.5cm}
\pgfputat{\pgfxy(-0.5,12)}{\pgfbox[center,center]{(i)}}
{\pgfsetstartarrow{\pgfarrowbar}
\pgfsetendarrow{\pgfarrowbar}
\pgfxyline(0,12)(2.2,12)
\pgfputat{\pgflabel{.5}{\pgfxy(0,12)}{\pgfxy(2.2,12)}{5pt}}{\pgfbox[center,base]{$V$}}
\pgfxyline(2.2,12)(6,12)
\pgfputat{\pgflabel{.5}{\pgfxy(3,12)}{\pgfxy(6,12)}{5pt}}{\pgfbox[center,base]{$T_{B,D}$}}
\pgfxyline(0,11)(2.5,11)
\pgfputat{\pgflabel{.5}{\pgfxy(0,11)}{\pgfxy(2.5,11)}{5pt}}{\pgfbox[center,base]{$T_{B,L}$}}
\pgfxyline(2.5,11)(5,11)
\pgfputat{\pgflabel{.5}{\pgfxy(2.5,11)}{\pgfxy(5,11)}{5pt}}{\pgfbox[center,base]{$T_{B,I}$}}
\pgfxyline(2.5,10.25)(6,10.25)
\pgfputat{\pgflabel{.5}{\pgfxy(2.5,10.25)}{\pgfxy(6,10.25)}{5pt}}{\pgfbox[center,base]{$W$}}
}

\pgfputat{\pgfxy(7,12)}{\pgfbox[center,center]{(ii)}}
{\pgfsetstartarrow{\pgfarrowbar}
\pgfsetendarrow{\pgfarrowbar}
\pgfxyline(7.5,12)(9.5,12)
\pgfputat{\pgflabel{.5}{\pgfxy(7.5,12)}{\pgfxy(9.5,12)}{5pt}}{\pgfbox[center,base]{$V$}}
\pgfxyline(9.5,12)(11.5,12)
\pgfputat{\pgflabel{.5}{\pgfxy(9.5,12)}{\pgfxy(11.5,12)}{5pt}}{\pgfbox[center,base]{$T_{B,D}$}}
\pgfxyline(7.5,11)(10,11)
\pgfputat{\pgflabel{.5}{\pgfxy(7.5,11)}{\pgfxy(10,11)}{5pt}}{\pgfbox[center,base]{$T_{B,L}$}}
\pgfxyline(10,11)(12.5,11)
\pgfputat{\pgflabel{.5}{\pgfxy(10,11)}{\pgfxy(12.5,11)}{5pt}}{\pgfbox[center,base]{$T_{B,I}$}}
\pgfxyline(10,10.25)(11.5,10.25)
\pgfputat{\pgflabel{.5}{\pgfxy(10,10.25)}{\pgfxy(11.5,10.25)}{5pt}}{\pgfbox[center,base]{$W$}}
}

\pgfputat{\pgfxy(-0.5,9)}{\pgfbox[center,center]{(iii)}}
{\pgfsetstartarrow{\pgfarrowbar}
\pgfsetendarrow{\pgfarrowbar}
\pgfxyline(0,9)(1,9)
\pgfputat{\pgflabel{.5}{\pgfxy(0,9)}{\pgfxy(1,9)}{5pt}}{\pgfbox[center,base]{$V$}}
\pgfxyline(1,9)(2.25,9)
\pgfputat{\pgflabel{.5}{\pgfxy(1,9)}{\pgfxy(2.25,9)}{5pt}}{\pgfbox[center,base]{$T_{B,D}$}}
\pgfxyline(0,8)(2.5,8)
\pgfputat{\pgflabel{.5}{\pgfxy(0,8)}{\pgfxy(2.5,8)}{5pt}}{\pgfbox[center,base]{$T_{B,L}$}}
\pgfxyline(2.5,8)(5,8)
\pgfputat{\pgflabel{.5}{\pgfxy(2.5,8)}{\pgfxy(5,8)}{5pt}}{\pgfbox[center,base]{$T_{B,I}$}}}
\end{pgfpicture}
\caption{Examples showing how $V$, $T_{B,D}$, $T_{B,L}$, $T_{B,I}$ and $W$ relate. In (i): $B$ dies naturally and is untraced; $B$'s active lifetime length is $T_{B,I}$. In (ii): $B$ is traced while active and has an active lifetime of length $W$. In (iii): $B$ is traced during its latent period and so has no offspring, $W=0$.}
\label{fig:TVWexamples}
\end{figure}

To calculate $m_{NU}$ and $m_{NN}$, consider a typical unnamed individual $A$ and a typical named offspring, $B$ say, of $A$. Let $V$ be the length of time between $B$'s birth and $A$'s death (which is when $B$ is named), so $V\sim\mbox{U}(0,\iota)$, and let $T_{B,D}$ and $T_{B,L}$ be the tracing delay and latent period of $B$, respectively. Hence the length of time after $B$'s birth that $B$ can be traced is $V+T_{B,D}$, but $B$ is not active (and so cannot have offspring) for the first $T_{B,L}$ time units, after which $B$ is active until the end of its lifetime or they are traced, whichever comes first. Letting $T_{D-L}=T_{B,D}-T_{B,L}$, $B$'s active lifetime is therefore $\min\{\iota,\max\{0,V+T_{D-L}\}\}$. See Figure \ref{fig:TVWexamples} for some examples; $T_{B,I}$ is $B$'s natural active lifetime length (here $\equiv\iota$) and $W$ (used in Section \ref{sec:expcalc}) is given by $\max\{0,V+T_{D-L}\}$.

Let $p_N$ be the probability that a named individual dies before the end of its tracing delay and so is untraced (as in example (i) in Figure \ref{fig:TVWexamples}). Then
\begin{align*}
p_N & = \mbox{P}\left(\iota\leq V+T_{D-L}\right)\\
    & = \mbox{P}\left(T_{D-L}>\iota\right)+\mbox{P}\left(\iota\leq V+T_{D-L},T_{D-L}\leq\iota\right)\\
    & = \mbox{P}\left(T_{D-L}>\iota\right)+\frac{1}{\iota}\mbox{E}\left[T_{D-L}1_{\{0\leq T_{D-L}\leq\iota\}}\right],
\end{align*}
since $\iota-V\sim\mbox{U}(0,\iota)$. (For an event, $F$ say, $1_F$ denotes its indicator function; i.e.~$1_F=1$ if the event $F$ occurs, and $1_F=0$ otherwise.) Such an individual behaves in the same manner as an unnamed individual. Thus, 
\begin{align*}
m_{NN} & = \lambda\pi_Rp\iota p_N,
\end{align*}
since traced individuals cannot name their offspring, and named but untraced individuals contribute $\lambda(1-\pi_Rp)\iota p_N$ to $m_{NU}$.

The remaining contribution to $m_{NU}$ comes from traced individuals (as in example (ii) in Figure \ref{fig:TVWexamples}) and is
\begin{align*}
&\quad \lambda\mbox{E}\left[(V+T_{D-L})1_{\left\{0<V+T_{D-L}<\iota\right\}}\right] \displaybreak[0]\\
& = \lambda\mbox{E}\left[(V+T_{D-L})1_{\left\{0<V+T_{D-L}<\iota,-\iota<T_{D-L}<0\right\}}\right]\\ 
&\quad+\lambda\mbox{E}\left[(V+T_{D-L})1_{\left\{0<V+T_{D-L}<\iota,0\leq T_{D-L}<\iota\right\}}\right] \displaybreak[0]\\
& = \lambda\mbox{E}\left[(\iota+T_{D-L}-(\iota-V))1_{\left\{0<\iota-V<\iota+T_{D-L},-\iota<T_{D-L}<0\right\}}\right]\\ 
&\quad+\lambda\mbox{E}\left[(\iota+T_{D-L}-(\iota-V))1_{\left\{T_{D-L}<\iota-V<\iota,0\leq T_{D-L}<\iota\right\}}\right] \displaybreak[0]\\
& = \lambda\mbox{E}_{T_{D-L}}\left[\mbox{E}_V\left[\left.(\iota+T_{D-L}-(\iota-V))1_{\left\{0<\iota-V<\iota+T_{D-L}\right\}}\right|T_{D-L}\right]1_{\left\{-\iota<T_{D-L}<0\right\}}\right]\\ 
&\quad+\lambda\mbox{E}_{T_{D-L}}\left[\mbox{E}_V\left[\left.(\iota+T_{D-L}-(\iota-V))1_{\left\{T_{D-L}<\iota-V<\iota\right\}}\right|T_{D-L}\right]1_{\left\{0\leq T_{D-L}<\iota\right\}}\right] \displaybreak[0]\\
& = \frac{\lambda}{2\iota}\mbox{E}\left[\left(\iota+T_{D-L}\right)^21_{\left\{-\iota<T_{D-L}<0\right\}}\right]\\
&\quad +\frac{\lambda}{2\iota}\left\{\iota^2\mbox{P}\left(0\leq T_{D-L}<\iota\right)-\mbox{E}\left[T_{D-L}^21_{\left\{0\leq T_{D-L}<\iota\right\}}\right]\right\},
\end{align*}
since $\iota-V\sim U(0,\iota)$ and, if $x\in(-\iota,0)$ and $U\sim\mbox{U}(0,\iota)$ then
\begin{align*}
\mbox{E}\left[(\iota+x-U)1_{\{0<U<\iota+x\}}\right] & = \frac{1}{\iota}\int_0^{\iota+x}(\iota+x-u)\hspace{1mm}du=\frac{(\iota+x)^2}{2\iota},
\end{align*}
while, if $x\in(0,\iota)$ and $U\sim\mbox{U}(0,\iota)$ then
\begin{align*}
\mbox{E}\left[(\iota+x-U)1_{\{x<U<\iota\}}\right] & = \frac{1}{\iota}\int_x^{\iota}(\iota+x-u)\hspace{1mm}du=\frac{1}{2\iota}(\iota^2-x^2). 
\end{align*}

Putting this together,
\begin{align*}
m_{NU} & = \lambda(1-\pi_Rp)\iota p_N+\frac{\lambda}{2\iota}\mbox{E}\left[\left(\iota+T_{D-L}\right)^21_{\left\{-\iota<T_{D-L}<0\right\}}\right]\\
&\quad +\frac{\lambda}{2\iota}\left\{\iota^2\mbox{P}\left(0\leq T_{D-L}<\iota\right)-\mbox{E}\left[T_{D-L}^21_{\left\{0\leq T_{D-L}<\iota\right\}}\right]\right\}.
\end{align*}

\subsection{Calculation of probability generating functions}
\label{subsec:pgfs}

Recall that the number of offspring of an unnamed individual has a Poisson distribution with mean $\lambda\iota$, and with probability $(1-\pi_R)$ these are all unnamed otherwise they are named independently with probability $p$. Thus,
\begin{align*}
f_U\left(s_U,s_N\right) & = \sum_{k=0}^\infty\frac{(\lambda\iota)^k}{k!}e^{-\lambda\iota}\left[\left(1-\pi_R\right)s_U^k+\pi_R\sum_{j=0}^k\binom{k}{j}p^j(1-p)^{k-j}s_U^{k-j}s_N^j\right]\\
						& = \sum_{k=0}^\infty\frac{(\lambda\iota)^k}{k!}e^{-\lambda\iota}\left[\left(1-\pi_R\right)s_U^k+\pi_R\left((1-p)s_U+ps_N\right)^k\right]\\
						& = \left(1-\pi_R\right)e^{-\lambda\iota\left(1-s_U\right)}+\pi_Re^{-\lambda\iota\left(1-(1-p)s_U-ps_N\right)}.
\end{align*}

To calculate $f_N\left(s_U,s_N\right)$, we consider the contributions of (a) named individuals who are untraced, (b) individuals traced during their latent period and (c) individuals traced after their latent period, corresponding to examples (i), (iii) and (ii), respectively, in Figure \ref{fig:TVWexamples}.

Recall that named individuals are untraced with probability $p_N$, with such individuals behaving like unnamed individuals, so their contribution to   $f_N\left(s_U,s_N\right)$ is $p_Nf_U\left(s_U,s_N\right)$.

Let $p_T$ be the probability a named individual is traced before the end of their latent period. Such an individual has no offspring (and so contributes $p_T$ to $f_N\left(s_U,s_N\right)$) and
\begin{align*}
p_T & = \mbox{P}\left(V+T_{D-L}\leq 0\right)\\
    & = \mbox{P}\left(T_{D-L}\leq -\iota\right)+\mbox{P}\left(V+T_{D-L}\leq 0,-\iota<T_{D-L}\leq 0 \right)\\
    & = \mbox{P}\left(T_{D-L}\leq -\iota\right)+\frac{1}{\iota}\mbox{E}\left[-T_{D-L}1_{\left\{-\iota<T_{D-L}\leq 0\right\}}\right].
\end{align*} 

The remaining contribution to $f_N\left(s_U,s_N\right)$ comes from individuals traced after their latent period and is
\begin{align*}
&\quad \mbox{E}\left[e^{-\lambda\left(V+T_{D-L}\right)\left(1-s_U\right)}1_{\left\{0<V+T_{D-L}<\iota\right\}}\right] \displaybreak[0]\\
& = \mbox{E}\left[e^{-\lambda\left(V+T_{D-L}\right)\left(1-s_U\right)}1_{\left\{0<V+T_{D-L}<\iota,-\iota<T_{D-L}<0\right\}}\right]\\ &\quad+\mbox{E}\left[e^{-\lambda\left(V+T_{D-L}\right)\left(1-s_U\right)}1_{\left\{0<V+T_{D-L}<\iota,0\leq T_{D-L}<\iota\right\}}\right] \displaybreak[0]\\
& = \lambda\mbox{E}\left[e^{-\lambda\left(\iota+T_{D-L}-(\iota-V)\right)\left(1-s_U\right)}1_{\left\{0<\iota-V<\iota+T_{D-L},-\iota<T_{D-L}<0\right\}}\right]\\ 
&\quad+\lambda\mbox{E}\left[e^{-\lambda\left(\iota+T_{D-L}-(\iota-V)\right)\left(1-s_U\right)}1_{\left\{T_{D-L}<\iota-V<\iota,0\leq T_{D-L}<\iota\right\}}\right] \displaybreak[0]\\
& = \lambda\mbox{E}_{T_{D-L}}\left[\mbox{E}_V\left[\left.e^{-\lambda\left(\iota+T_{D-L}-(\iota-V)\right)\left(1-s_U\right)}1_{\left\{0<\iota-V<\iota+T_{D-L}\right\}}\right|T_{D-L}\right]1_{\left\{-\iota<T_{D-L}<0\right\}}\right]\\ 
&\quad+\lambda\mbox{E}_{T_{D-L}}\left[\mbox{E}_V\left[\left.e^{-\lambda\left(\iota+T_{D-L}-(\iota-V)\right)\left(1-s_U\right)}1_{\left\{T_{D-L}<\iota-V<\iota\right\}}\right|T_{D-L}\right]1_{\left\{0\leq T_{D-L}<\iota\right\}}\right] \displaybreak[0]\\
& = \frac{1}{\lambda\iota}\left\{\mbox{P}\left(-\iota<T_{D-L}<0\right)-\mbox{E}\left[e^{-\lambda\left(\iota+T_{D-L}\right)\left(1-s_U\right)}1_{\left\{-\iota<T_{D-L}<0\right\}}\right]\right\}\\
&\quad +\frac{1}{\lambda\iota}\left\{\mbox{E}\left[e^{-\lambda T_{D-L}\left(1-s_U\right)}1_{\left\{0\leq T_{D-L}<\iota\right\}}\right]-e^{-\lambda\iota\left(1-s_U\right)}\mbox{P}\left(0\leq T_{D-L}<\iota\right)\right\},
\end{align*}
since $\iota-V\sim U(0,\iota)$ and, if $x\in(-\iota,0)$ and $U\sim\mbox{U}(0,\iota)$ then, for $\theta\geq0$,
\begin{align*}
\mbox{E}\left[e^{-\theta(\iota+x-U)}1_{\{0<U<\iota+x\}}\right] & = \frac{1}{\iota}\int_0^{\iota+x}e^{-\theta(\iota+x-u)}\hspace{1mm}du=\frac{1}{\theta\iota}\left(1-e^{-\theta(\iota+x)}\right),
\end{align*}
while, if $x\in(0,\iota)$ and $U\sim\mbox{U}(0,\iota)$ then, for $\theta\geq0$,
\begin{align*}
\mbox{E}\left[e^{-\theta(\iota+x-U)}1_{\{x<U<\iota\}}\right] & = \frac{1}{\iota}\int_x^{\iota}e^{-\theta(\iota+x-u)}\hspace{1mm}du=\frac{1}{\theta\iota}\left(e^{-\theta x}-e^{-\theta\iota}\right). 
\end{align*}

Putting this all together,
\begin{align*}
f_N\left(s_U,s_N\right) & = p_T +p_Nf_U\left(q_U,q_N\right)\\
&\quad + \frac{1}{\lambda\iota}\left\{\mbox{P}\left(-\iota<T_{D-L}<0\right)-\mbox{E}\left[e^{-\lambda\left(\iota+T_{D-L}\right)\left(1-s_U\right)}1_{\left\{-\iota<T_{D-L}<0\right\}}\right]\right.\\
&\quad +\left.\mbox{E}\left[e^{-\lambda T_{D-L}\left(1-s_U\right)}1_{\left\{0\leq T_{D-L}<\iota\right\}}\right]-e^{-\lambda\iota\left(1-s_U\right)}\mbox{P}\left(0\leq T_{D-L}<\iota\right)\right\}.
\end{align*}

\section{Embedded Galton-Watson process (EGWP)}
\label{sec:EGWP}

In Section \ref{sec:const}, we saw that in the special case where $T_I\equiv\iota$ and $\pi_T=0$, named individuals are independent, so we can analyse the threshold behaviour of the MBDP (modified birth-death process) via a two-type Galton-Watson process. This approach is not possible in general since named individuals in the MBDP with the same parent are co-dependent on the length of their parent's active lifetime. Additionally the behaviour of named individuals is determined by how far down a naming chain they are.

However, in the MBDP, unnamed individuals are independent of one another and have active lifetimes of their full natural length. The threshold behaviour of the modified birth-death process can be obtained then by considering the embedded single-type discrete-time Galton-Watson process (EGWP) describing unnamed individuals, in which the offspring of a given individual are (a) their immediate unnamed offspring and (b) unnamed descendants who are separated from the given individual, in the family tree, only by named individuals.

\begin{figure}[h!]
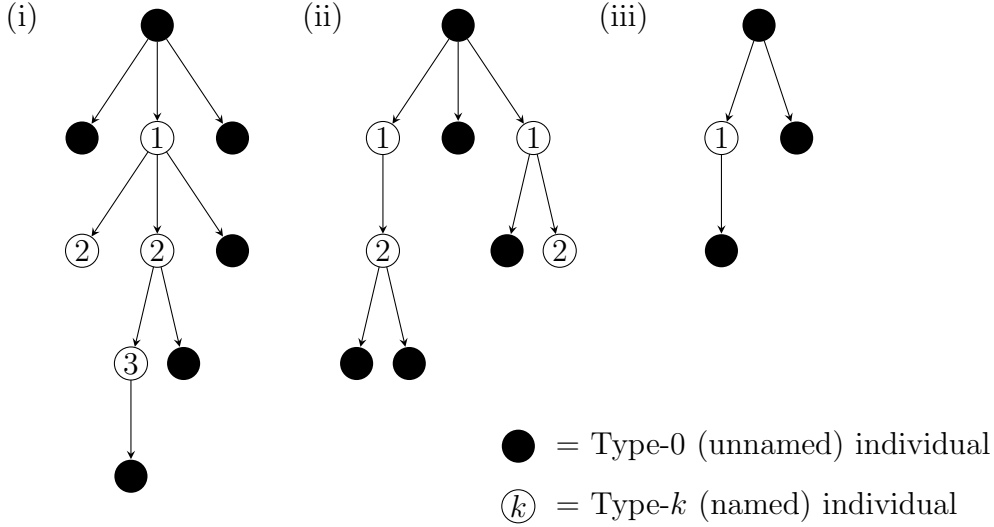

\begin{pgfpicture}{-0.3cm}{3cm}{12cm}{10.5cm}
\pgfputat{\pgfxy(0.2,9.9)}{\pgfbox[center,center]{(i)}}
\pgfnodecircle{Node9}[fill]{\pgfxy(2,9.8)}{0.22cm}
\pgfnodecircle{Node10}[stroke]{\pgfrelative{\pgfxy(0,-1.5)}{\pgfnodecenter{Node9}}}{0.22cm}
\pgfputat{\pgfrelative{\pgfxy(0,0)}{\pgfnodecenter{Node10}}}{\pgfbox[center,center]{$1$}}
\pgfsetendarrow{\pgfarrowsingle}
\pgfnodesetsepstart{0pt}
\pgfnodesetsepend{1pt}
\pgfnodeconnline{Node9}{Node10}
\pgfnodecircle{Node11}[fill]{\pgfrelative{\pgfxy(-1,-1.5)}{\pgfnodecenter{Node9}}}{0.22cm}
\pgfsetendarrow{\pgfarrowsingle}
\pgfnodesetsepstart{0pt}
\pgfnodesetsepend{1pt}
\pgfnodeconnline{Node9}{Node11}
\pgfnodecircle{Node12}[fill]{\pgfrelative{\pgfxy(1,-1.5)}{\pgfnodecenter{Node9}}}{0.22cm}
\pgfsetendarrow{\pgfarrowsingle}
\pgfnodesetsepstart{0pt}
\pgfnodesetsepend{1pt}
\pgfnodeconnline{Node9}{Node12}
\pgfnodecircle{Node13}[stroke]{\pgfrelative{\pgfxy(-1,-1.5)}{\pgfnodecenter{Node10}}}{0.22cm}
\pgfputat{\pgfrelative{\pgfxy(0,0)}{\pgfnodecenter{Node13}}}{\pgfbox[center,center]{2}}
\pgfsetendarrow{\pgfarrowsingle}
\pgfnodesetsepstart{0pt}
\pgfnodesetsepend{1pt}
\pgfnodeconnline{Node10}{Node13}
\pgfnodecircle{Node14}[stroke]{\pgfrelative{\pgfxy(0,-1.5)}{\pgfnodecenter{Node10}}}{0.22cm}
\pgfputat{\pgfrelative{\pgfxy(0,0)}{\pgfnodecenter{Node14}}}{\pgfbox[center,center]{2}}
\pgfsetendarrow{\pgfarrowsingle}
\pgfnodesetsepstart{0pt}
\pgfnodesetsepend{1pt}
\pgfnodeconnline{Node10}{Node14}
\pgfnodecircle{Node18}[fill]{\pgfrelative{\pgfxy(1,-1.5)}{\pgfnodecenter{Node10}}}{0.22cm}
\pgfsetendarrow{\pgfarrowsingle}
\pgfnodesetsepstart{0pt}
\pgfnodesetsepend{1pt}
\pgfnodeconnline{Node10}{Node18}
\pgfnodecircle{Node15}[stroke]{\pgfrelative{\pgfxy(-0.35,-1.5)}{\pgfnodecenter{Node14}}}{0.22cm}
\pgfputat{\pgfrelative{\pgfxy(0,0)}{\pgfnodecenter{Node15}}}{\pgfbox[center,center]{3}}
\pgfsetendarrow{\pgfarrowsingle}
\pgfnodesetsepstart{0pt}
\pgfnodesetsepend{1pt}
\pgfnodeconnline{Node14}{Node15}
\pgfnodecircle{Node16}[fill]{\pgfrelative{\pgfxy(0.35,-1.5)}{\pgfnodecenter{Node14}}}{0.22cm}
\pgfsetendarrow{\pgfarrowsingle}
\pgfnodesetsepstart{0pt}
\pgfnodesetsepend{1pt}
\pgfnodeconnline{Node14}{Node16}
\pgfnodecircle{Node17}[fill]{\pgfrelative{\pgfxy(0,-1.5)}{\pgfnodecenter{Node15}}}{0.22cm}
\pgfsetendarrow{\pgfarrowsingle}
\pgfnodesetsepstart{0pt}
\pgfnodesetsepend{1pt}
\pgfnodeconnline{Node15}{Node17}

\pgfputat{\pgfxy(4.2,9.9)}{\pgfbox[center,center]{(ii)}}
\pgfnodecircle{Node1}[fill]{\pgfxy(6,9.8)}{0.22cm}
\pgfnodecircle{Node2}[fill]{\pgfrelative{\pgfxy(0,-1.5)}{\pgfnodecenter{Node1}}}{0.22cm}
\pgfsetendarrow{\pgfarrowsingle}
\pgfnodesetsepstart{0pt}
\pgfnodesetsepend{1pt}
\pgfnodeconnline{Node1}{Node2}
\pgfnodecircle{Node3}[stroke]{\pgfrelative{\pgfxy(-1,-1.5)}{\pgfnodecenter{Node1}}}{0.22cm}
\pgfputat{\pgfrelative{\pgfxy(0,0)}{\pgfnodecenter{Node3}}}{\pgfbox[center,center]{1}}
\pgfsetendarrow{\pgfarrowsingle}
\pgfnodesetsepstart{0pt}
\pgfnodesetsepend{1pt}
\pgfnodeconnline{Node1}{Node3}
\pgfnodecircle{Node4}[stroke]{\pgfrelative{\pgfxy(1,-1.5)}{\pgfnodecenter{Node1}}}{0.22cm}
\pgfputat{\pgfrelative{\pgfxy(0,0)}{\pgfnodecenter{Node4}}}{\pgfbox[center,center]{1}}
\pgfsetendarrow{\pgfarrowsingle}
\pgfnodesetsepstart{0pt}
\pgfnodesetsepend{1pt}
\pgfnodeconnline{Node1}{Node4}
\pgfnodecircle{Node5}[stroke]{\pgfrelative{\pgfxy(0,-1.5)}{\pgfnodecenter{Node3}}}{0.22cm}
\pgfputat{\pgfrelative{\pgfxy(0,0)}{\pgfnodecenter{Node5}}}{\pgfbox[center,center]{2}}
\pgfsetendarrow{\pgfarrowsingle}
\pgfnodesetsepstart{0pt}
\pgfnodesetsepend{1pt}
\pgfnodeconnline{Node3}{Node5}
\pgfnodecircle{Node6}[fill]{\pgfrelative{\pgfxy(-0.35,-1.5)}{\pgfnodecenter{Node5}}}{0.22cm}
\pgfsetendarrow{\pgfarrowsingle}
\pgfnodesetsepstart{0pt}
\pgfnodesetsepend{1pt}
\pgfnodeconnline{Node5}{Node6}
\pgfnodecircle{Node7}[fill]{\pgfrelative{\pgfxy(0.35,-1.5)}{\pgfnodecenter{Node5}}}{0.22cm}
\pgfsetendarrow{\pgfarrowsingle}
\pgfnodesetsepstart{0pt}
\pgfnodesetsepend{1pt}
\pgfnodeconnline{Node5}{Node7}
\pgfnodecircle{Node8}[fill]{\pgfrelative{\pgfxy(-0.35,-1.5)}{\pgfnodecenter{Node4}}}{0.22cm}
\pgfsetendarrow{\pgfarrowsingle}
\pgfnodesetsepstart{0pt}
\pgfnodesetsepend{1pt}
\pgfnodeconnline{Node4}{Node8}
\pgfnodecircle{Node9a}[stroke]{\pgfrelative{\pgfxy(0.35,-1.5)}{\pgfnodecenter{Node4}}}{0.22cm}
\pgfputat{\pgfrelative{\pgfxy(0,0)}{\pgfnodecenter{Node9a}}}{\pgfbox[center,center]{2}}
\pgfsetendarrow{\pgfarrowsingle}
\pgfnodesetsepstart{0pt}
\pgfnodesetsepend{1pt}
\pgfnodeconnline{Node4}{Node9a}

\pgfputat{\pgfxy(8.2,9.9)}{\pgfbox[center,center]{(iii)}}
\pgfnodecircle{Node19}[fill]{\pgfxy(10,9.8)}{0.22cm}
\pgfnodecircle{Node20}[stroke]{\pgfrelative{\pgfxy(-0.5,-1.5)}{\pgfnodecenter{Node19}}}{0.22cm}
\pgfputat{\pgfrelative{\pgfxy(0,0)}{\pgfnodecenter{Node20}}}{\pgfbox[center,center]{1}}
\pgfsetendarrow{\pgfarrowsingle}
\pgfnodesetsepstart{0pt}
\pgfnodesetsepend{1pt}
\pgfnodeconnline{Node19}{Node20}
\pgfnodecircle{Node21}[fill]{\pgfrelative{\pgfxy(0.5,-1.5)}{\pgfnodecenter{Node19}}}{0.22cm}
\pgfsetendarrow{\pgfarrowsingle}
\pgfnodesetsepstart{0pt}
\pgfnodesetsepend{1pt}
\pgfnodeconnline{Node19}{Node21}
\pgfnodecircle{Node22}[fill]{\pgfrelative{\pgfxy(0,-1.5)}{\pgfnodecenter{Node20}}}{0.22cm}
\pgfsetendarrow{\pgfarrowsingle}
\pgfnodesetsepstart{0pt}
\pgfnodesetsepend{1pt}
\pgfnodeconnline{Node20}{Node22}

\pgfnodecircle{Node23}[fill]{\pgfxy(6.8,4.2)}{0.22cm}
\pgfnodecircle{Node24}[stroke]{\pgfrelative{\pgfxy(0,-0.8)}{\pgfnodecenter{Node23}}}{0.22cm}
\pgfputat{\pgfrelative{\pgfxy(0,0)}{\pgfnodecenter{Node24}}}{\pgfbox[center,center]{$k$}}
\pgfputat{\pgfrelative{\pgfxy(0.5,0)}{\pgfnodecenter{Node23}}}{\pgfbox[left,center]{= Type-$0$ (unnamed) individual}}
\pgfputat{\pgfrelative{\pgfxy(0.5,0)}{\pgfnodecenter{Node24}}}{\pgfbox[left,center]{= Type-$k$ (named) individual}}
\end{pgfpicture}
\caption[Some example realizations of $R$.]{Some example realizations of $R$.\\
In (i) $R=5$, $R_1=1$, $R_2=1$, $R_3=1$, $R^{(0)}=2$, $R^{(1)}=3$, $R^{(2)}=4$, $R^{(3)}=5$.\\
In (ii) $R=4$, $R_1=1$, $R_2=2$, $R^{(0)}=1$, $R^{(1)}=2$, $R^{(2)}=4$.\\
In (iii) $R=2$, $R_1=1$, $R^{(0)}=1$, $R^{(1)}=2$.}
\label{fig:Rexamples}
\end{figure}

To analyse the MBDP we focus attention on the offspring random variable, $R$ say, in the EGWP, by obtaining expressions for its mean, $R_U=\ez[R]$, which we call a type-reproduction number, following Heesterbeek \& Roberts \cite{heesterbeek-roberts}. For some example realizations of $R$, see Figure \ref{fig:Rexamples} (notation such as $R_1$, $R^{(0)}$ etc is defined below). Standard results from branching process theory (e.g.~Haccou {\em et al.}~\cite{haccou}, Theorem 5.2) tell us that the EGWP will die out with probability $1$ if and only if $R_U\leq1$, and that if $R_U>1$ the extinction probability of the EGWP is $p_E^m$ where $p_E$ is the solution of $s=H(s)$ in $(0,1)$, where $H(s)=\mbox{E}\left[s^R\right]$.

Let a named individual who is separated from an ancestor who was unnamed in the family tree of the MBDP by $k-1$ named individuals only be called a type-$k$ individual. Hence, type-$1$ individuals are the named immediate offspring of unnamed individuals, type-$2$ individuals are the named immediate offspring of the named immediate offspring of unnamed individuals, and so on. Type-$0$ individuals are unnamed individuals. We refer to the type-$k$ individuals who are separated from a given unnamed individual in the family tree only by named individuals as that unnamed individual's immediate type-$k$ descendants.

We obtain $R_U$ by considering again a typical unnamed individual, $A$, in the MBDP. Let $R_i$ be the total number of unnamed immediate offspring of all the immediate type-$i$ ($i=1,2,\ldots$) descendants of $A$ and let $R_{U,i}=\mbox{E}\left[R_i\right]$. Let $R^{(k)}$ be the total number of unnamed immediate offspring of $A$ and all of its immediate descendants of up to and including type-$k$. Then for $k=1,2,\ldots$,
\begin{align}
R^{(k)} & =  R^{(0)} + \sum_{i=1}^kR_i,
\label{eqn:Rk}
\end{align}
where $R^{(0)}$ is the number of unnamed immediate offspring of $A$. Examples of $R_i$ and $R^{(k)}$ are given in Figure \ref{fig:Rexamples}. We have that $R^{(k)}\uparrow R$ as $k\to\infty$, so by the monotone convergence theorem,
\begin{align}
R_U & = \lim_{k\to\infty}R_U^{(k)},
\label{eqn:RU}
\end{align}
where $R_U^{(k)}=\ez\left[R^{(k)}\right]$.

In the MBDP, the unnamed individual $A$ produces a cluster of named individuals (which may be of size zero), and $A$'s immediate type-$k$ descendants are the $k$th generation individuals within this naming cluster. The offspring of $A$ in the EGWP ($R$) is all of the unnamed offspring of $A$ and of the individuals in the naming cluster, with $R_i$ being the total number of unnamed offspring of the $i$th generation.

We revisit the case where $T_I\equiv\iota$ and $\pi_T=0$. Note that the two-type Galton-Watson process goes extinct if and only if the EGWP goes extinct, so the extinction probability is the same (i.e.~$p_E=q_U$) and $R_0$ and $R_U$ have the same threshold (see Roberts and Heesterbeek \cite{roberts-heesterbeek}). 

An unnamed individual is expected to have $m_{UN}$ immediate named offspring, each of whom is expected to have $m_{NN}$ immediate named offspring, and so on. Hence the expected number of (named) immediate type-$i$ descendants ($i=1,\ldots$) of an unnamed individual is $m_{UN}m_{NN}^{i-1}$. Since unnamed and named individuals are expected to have $m_{UU}$ and $m_{NU}$ unnamed offspring, respectively, $\mbox{E}\left[R^{(0)}\right]=m_{UU}$, and, for $i=1,2,\ldots,$ $\mbox{E}\left[R_i\right]=m_{UN}m_{NN}^{i-1}m_{NU}$. Therefore $R_U^{(k)}=m_{UU}+\sum_{i=1}^km_{UN}m_{NN}^{i-1}m_{NU}$, whence, using Eqn.~(\ref{eqn:RU}),
\begin{align*}
R_U & = \left\{\begin{array}{ll}
    		m_{UU}+\frac{m_{UN}m_{NU}}{1-m_{NN}} & \mbox{if $m_{NN}<1$,}\\
    		\infty                               & \mbox{if $m_{NN}\geq1$.}\end{array}\right.
\end{align*}

We assume generally that siblings experience independent delays. Note that if we were to assume instead that siblings experience delays of the same length, then this will affect offspring distributions, but offspring means will remain the same. Hence $R_U$ would be the same but the extinction probability would be different. Since $R_U$ is unchanged and $R_U$ and $R_0$ have the same threshold, $R_0$ (as given by Eqn.~(\ref{eqn:R0})) is still a threshold parameter when siblings experience the same delay.

\section{Exponential infectious period}
\label{sec:exp}

In this section we assume that $T_I\sim\mbox{Exp}(\gamma)$, i.e.~the infectious period of each individual in the epidemic model is exponentially distributed with mean $\gamma^{-1}$.

\subsection{Calculation of threshold parameter}
\label{sec:expcalc}

Recall that in the constant infectious period case (with $\pi_T=0$) that $R_U<\infty$ if and only if $m_{NN}=\lambda\pi_Rp\iota p_N<1$, thus, if all other parameters are held fixed, $R_U<\infty$ if and only if $\lambda<\lambda^*=\left(\pi_Rp\iota p_N\right)^{-1}$. 

The following lemma is a similar result that holds in the current setting, and is proved in \ref{app:RUinfty}.
\begin{lem}
If $\mbox{P}\left(T_D>T_L\right)=0$, then $R_U<\infty$ for all $\lambda<\infty$. If $\mbox{P}\left(T_D>T_L\right)>0$, then there exists some $\lambda^*>\gamma$ such that, if $\lambda<\lambda^*$, $R_U<\infty$, while if $\lambda>\lambda^*$, $R_U=\infty$.
\end{lem}

We now assume the delays are exponentially distributed with mean $\xi^{-1}$, i.e.~$T_D\sim\mbox{Exp}(\xi)$. We let the latent period have some arbitrary but specified distribution $T_L$, which for ease of exposition we assume to be continuous with probability distribution function $f_L(t)$ ($t>0$) and moment-generating function $\phi_L(\theta)=\mbox{E}\left[e^{-\theta T_L}\right]=\int_0^\infty e^{-\theta t}f_L(t)\hspace{1mm}dt$ ($\theta\geq0$). Results extend easily to the case where $T_L$ is discrete. With exponential delays we are now in the setting where $\mbox{P}\left(T_D>T_L\right)>0$.

Let $T_{A,I}\sim\mbox{Exp}(\gamma)$ denote the active lifetime of a typical unnamed individual, $A$, in the MBDP and, for $k=0,1,\ldots$, let
\begin{align}
h_k(t) & = \mbox{E}\left[R^{(k)}|T_{A,I}=t\right],
\label{eqn:hkt}
\end{align}
so
\begin{align*}
R_U & = \lim_{k\to\infty}\int_0^\infty\gamma e^{-\gamma t}h_k(t)\hspace{1mm}dt.
\end{align*}

Recall that $R^{(0)}$ is the number of unnamed immediate offspring of $A$. Hence, with probability $\pi_R$, $\left(R^{(0)}|T=t\right)\sim\mbox{Poisson}(\lambda(1-p)t)$, otherwise $\left(R^{(0)}|T=t\right)\sim\mbox{Poisson}(\lambda t)$, so
\begin{align}
h_0(t)=\mbox{E}\left[\left.R^{(0)}\right|T_{A,I}=t\right] & = \lambda t\left\{\pi_R(1-p)+(1-\pi_R)\right\}= \lambda\left(1-\pi_Rp\right)t.
\end{align}

Let $N_1$ denote the number of named immediate (i.e.~type-$1$) offspring of $A$ and, for $j=1,2,\ldots,N_1$, let $Z^{(k)}_j$ be the total number of (unnamed) descendants from the $i$th such (arbitrarily ordered) individual that contribute to $R^{(k)}$. Thus
\begin{align*}
\sum_{i=1}^kR_i & = \sum_{j=1}^{N_1}Z^{(k)}_j,
\end{align*}
where the sum on the right is zero if $N_1=0$. Now, with probability $1-\pi_R$, $N_1=0$, otherwise $\left(N_1|T_{A,I}=t\right)\sim\mbox{Poisson}(\lambda pt)$ and, conditional upon $N_1=n$ and $T_{A,I}=t$, the birth times of these $n$ individuals can be obtained by sampling $n$ independent $U(0,t)$ random variables, so the $Z^{(k)}_j$ are independent and identically distributed, with common distribution $Z^{(k)}$, say. Thus
\begin{align}
\mbox{E}\left[\sum_{i=1}^kR_i|T_{A,I}=t\right] & = \sum_{n=0}^{\infty}\mbox{P}\left(\left.N_1=n\right|T_{A,I}=t\right)n\mbox{E}\left[Z^{(k)}|T_{A,I}=t\right] \nonumber\\
                                         & = \lambda\pi_R ptg_k(t),
\label{eqn:sumRi}
\end{align}
where $g_k(t)=\mbox{E}\left[Z^{(k)}|T_{A,I}=t\right]$ ($k=1,2,\ldots$), so, using Eqns.~(\ref{eqn:Rk}) and (\ref{eqn:hkt})-(\ref{eqn:sumRi}), for $k=1,2,\ldots$
\begin{eqnarray}
h_k(t) = \lambda\left(1-\pi_Rp\right)t+\lambda\pi_Rptg_k(t) \hspace{1cm} (t>0).
\label{eqn:hk1}
\end{eqnarray}

Consider individual $B$, a typical named offspring of $A$, with natural active lifetime length $T_{B,I}$ and let $V$ be the length of time between $B$'s birth and $A$'s death (which is when $B$ is named), so
\begin{align*}
(V|T_{A,I}=t) & \sim\mbox{U}(0,t).
\end{align*}

Let $T_{B,D}$ and $T_{B,L}$ be the tracing delay and latent period of $B$, respectively. The length of time after $B$'s birth that $B$ can be traced is $V+T_{D,B}$. For the first $T_{B,L}$ time units, $B$ is not active (and is not having offspring), and after this they are active until the end of their natural active lifetime or they are traced, whichever comes first. Let $W=\max\left\{0,V+T_{B,D}-T_{B,L}\right\}$, so that $B$'s actual active lifetime is $\min\{T_{B,I},W\}$. Note that $W=0$ if $B$ is traced before coming active, otherwise $W$ represents the length of time between $B$ becoming active and the end of the delay in tracing $B$. See Figure \ref{fig:TVWexamples} for some examples.

Conditioning on $T_{A,I}=t$, $T_{B,D}=t_D$ and $T_{B,L}=t_L$, we consider two cases: (i) $t_D>t_L$ and (ii) $t_D\leq t_L$. In case (i), $W$ is uniformly distributed over $(t_D-t_L,t+t_D-t_L)$. In case (ii), $W$ has a probability mass of $\frac{t_L-t_D}{t}$ at $0$ with the remaining mass uniformly distributed over $(0,t+t_D-t_L)$.

If $T_{B,I}>W$ then $B$ is traced and so can name its immediate offspring with probability $\pi_T$, otherwise all $B$'s immediate offspring are unnamed, so
\begin{align*}
&\left(\left.Z^{(k)}\right|W=w,T_{B,I}>w\right) \eqd \left(R^{(k-1)}|T_{A,I}=w\right) \hspace{0.55cm} \mbox{if $B$ is interviewed},\\
&\mbox{E}\left[Z^{(k)}|W=w,T_{B,I}>w\right] =\lambda w \hspace{3.2cm} \mbox{if $B$ is not interviewed}.
\end{align*}
Conversely, if $T_{B,I}<W$ then $B$ dies before it can be traced as a result of being named by $A$, so $B$ dies naturally and can name its immediate offspring with probability $\pi_R$, otherwise all $B$'s immediate offspring are unnamed, so, if $t_I\leq w$,
\begin{align*}
&\left(\left.Z^{(k)}\right|W=w,T_{B,I}=t_I\right) \eqd \left(R^{(k-1)}|T_{A,I}=t_I\right) \hspace{0.55cm} \mbox{if $B$ is interviewed},\\
&\mbox{E}\left[Z^{(k)}|W=w,T_{B,I}=t_I\right] =\lambda t_I \hspace{3.2cm} \mbox{if $B$ is not interviewed}.
\end{align*}

Putting this all together yields
\begin{align}
\begin{split}
g_k(t) & = \int_{t_D=0}^\infty\xi e^{-\xi t_D}\int_{t_L=0}^{t_D} \int_{w=t_D-t_L}^{t+t_D-t_L}\frac{1}{t}\chi(w)f_L\left(t_L\right)\hspace{1mm}dw\hspace{1mm}dt_L\hspace{1mm}dt_D \\
       & \quad + \int_{t_D=0}^\infty\xi e^{-\xi t_D}\int_{t_L=t_D}^{t_D+t}\int_{w=0}^{t+t_D-t_L}\frac{1}{t}\chi(w)f_L\left(t_L\right)\hspace{1mm}dw \hspace{1mm}dt_L\hspace{1mm}dt_D,
\end{split}
\label{eqn:gk1}
\end{align}
where, for $w>0$,
\begin{align*}
\chi(w) = \mbox{E}\left[\left.Z^{(k)}\right|W=w\right] & = \lambda\left(1-\pi_T\right)we^{-\gamma w}+\pi_T e^{-\gamma w} h_{k-1}(w)\\
        & \quad +\int_{u=0}^w\gamma e^{-\gamma u}\left(\lambda\left(1-\pi_R\right)u+\pi_Rh_{k-1}(u)\right)du.
\end{align*}

Let $L_k(\theta) = \int_0^\infty e^{-\theta t}h_k(t)\hspace{1mm}dt$. Assume $\gamma=1$ without loss of generality. Then
\begin{align*}
R_U & = \lim_{k\to\infty}\int_0^\infty e^{-t}h_k(t)\hspace{1mm}dt\\
    & = \lim_{k\to\infty}L_k(1).
\end{align*}

It follows from $h_0(t)=\lambda\left(1-\pi_Rp\right)t$ that $L_0(\theta)=\lambda(1-\pi_Rp)/\theta^2$, and it follows (details in \ref{app:Lkthetadetails}) from Eqn.~(\ref{eqn:hk1}) that
\begin{align}
L_k(\theta)  & = \frac{\lambda(1-\pi_Rp)}{\theta^2}+ \frac{\lambda\pi_R p\xi}{\theta(\xi-\theta)} \left\{\phi_L(\theta)\psi(\theta)-\phi_L(\xi)\psi(\xi)\right\},
\label{eqn:Lktheta}
\end{align}
assuming that $\theta\neq\xi$ (the case where $\theta=\xi$ can be treated by taking the limit as $\xi\to\theta$), and where, for $\theta\geq0$,
\begin{align*}
\psi(\theta)& = \int_{0}^\infty e^{-\theta w}\chi(w) \hspace{1mm}dw = \frac{\lambda(1-\pi_T)}{(\theta+1)^2}+\frac{\lambda(1-\pi_R)}{\theta(\theta+1)^2}+\frac{\pi_T\theta+\pi_R}{\theta}L_{k-1}(\theta+1).
\end{align*}

Let $x_k(\theta)=L_k(\theta+1)$. Suppose that $\xi$ is non-integer. Then it follows from Eqn.~(\ref{eqn:Lktheta}) (details in \ref{app:xkdetails}) that
\begin{eqnarray}
x_k(\theta) = \sum_{i=0}^{k-1}c_i(\theta) \left(a_{i+1}(\theta)-\rho_{i+1}(\theta)x_{k-i-1}(\xi)\right)+\frac{c_k(\theta)\lambda(1-\pi_Rp)}{(k+1+\theta)^2},
\label{eqn:xk}
\end{eqnarray}
where, for $\theta\geq0$, $c_0(\theta)=1$ and, for $j=1,2,\ldots$,
\begin{align*}
a_j(\theta) & = \frac{\lambda(1-\pi_Rp)}{(j+\theta)^2}\\
& \quad + \frac{\lambda^2\pi_Rp\xi\left(1-\pi_T+(1-\pi_R)/(j+\theta)\right)}{(j+\theta)(\xi-j-\theta)}\left(\frac{\phi_L(j+\theta)}{(j+\theta+1)^2}-\frac{\phi_L(\xi)}{(\xi+1)^2}\right), \displaybreak[0]\\
\rho_j(\theta) & = \frac{\lambda\pi_R p(\pi_T\xi+\pi_R)\phi_L(\xi)}{(j+\theta)(\xi-j-\theta)}, \displaybreak[0]\\
b_j(\theta) & = \frac{\lambda\pi_R p\xi(\pi_T(j+\theta)+\pi_R)\phi_L(j+\theta)}{(j+\theta)^2(\xi-j-\theta)}, \displaybreak[0]\\
c_j(\theta) &= \prod_{i=1}^jb_i(\theta).
\end{align*}
The case where $\xi$ is an integer can be considered by the taking the limit as $\xi$ approaches that integer value in the above expressions.

Let $x(\theta)=\lim_{k\to\infty}x_k(\theta)$ and let $y=\lim_{k\to\infty}L_k(\xi+1)=x(\xi)$. It is shown in \ref{app:xthetadetails} that $y<\infty$ if and only if $R_U<\infty$.

Suppose that $\lambda<\lambda^*$ (so $R_U<\infty$). Letting $k\to\infty$ in Eqn.~(\ref{eqn:xk}) yields (details in \ref{app:xthetadetails}) that
\begin{eqnarray}
x(\theta) = \sum_{i=0}^{\infty}c_i(\theta) \left(a_{i+1}(\theta)-\rho_{i+1}(\theta)y\right).
\label{eqn:xtheta}
\end{eqnarray}
By setting $\theta=\xi$ and rearranging, we have that $y=y^*$, where
\begin{align}
y^* = \frac{\sum_{i=0}^{\infty}c_i(\xi)a_{i+1}(\xi)}{1+\sum_{i=0}^{\infty}c_i(\xi)\rho_{i+1}(\xi)}.
\label{eqn:ystar}
\end{align}

Noting that $R_U=x(0)$, setting $\theta=0$ in Eqn.~(\ref{eqn:xtheta}) gives $R_U$ for $\lambda<\lambda^*$, so
\begin{align}
R_U & = \left\{\begin{array}{ll}
\sum_{i=0}^{\infty}c_i(0) \left(a_{i+1}(0)-\rho_{i+1}(0)y^*\right) & \mbox{if $\lambda<\lambda^*$,}\\
\infty & \mbox{if $\lambda>\lambda^*$.}\end{array}\right.
\label{eqn:RUexp}
\end{align}

If we assume instead that siblings experience the same delay then the $Z_j^{(k)}$ are no longer independent. However, due to the linearity of expectations, Eqn.~(\ref{eqn:sumRi}) still holds, so $R_U$ is unchanged.

By conditioning on $A$'s lifetime and its number of immediate offspring, we can obtain an integral equation involving $\mbox{E}\left[\left.s^R\right|T_{A,I}=t\right]$, details of which may be found in Chapter 3 of Knock \cite{knock} for the case when $\pi_R=1$, $\pi_T=0$ and there is no latent period or delay (but can be readily extended to more general $(\pi_R,\pi_T)$ and incorporating latent periods and delays). However, it does not appear possible to solve this equation in order to obtain $H(s)$, and hence explicitly calculate the extinction probability. Unlike $R_U$, the extinction probability is affected by assuming instead that siblings experience the same delay.

\subsection{Illustration of convergence of $R_U$}

We have seen that $R_U$ is finite for $\lambda<\lambda^*$ and infinite for $\lambda>\lambda^*$. Here we present an example to illustrate this phenomenon and how one finds $\lambda^*$ in practice.

\begin{figure}[h!]
\centering
\centerline{\includegraphics[width=1.1\textwidth]{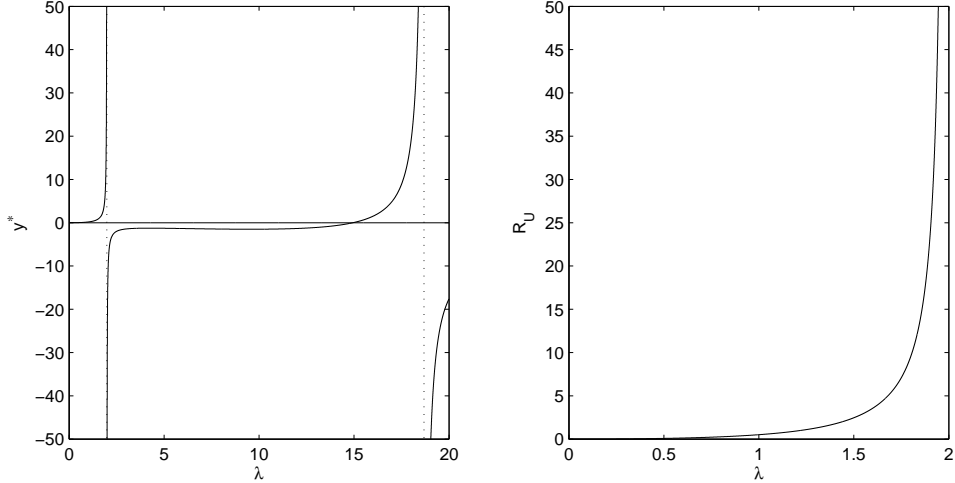}}
\caption{$y^*$ and $R_U$ varying with $\lambda$, when $T_I\sim\mbox{Exp}(1)$, $T_D\sim\mbox{Exp}(0.7)$, $T_L\equiv0$, $\pi_R=1$, $\pi_T=0$, and $p=1$.}
\label{fig:convergence}
\end{figure}

Figure \ref{fig:convergence} shows a plot of $y^*$ varying with $\lambda$ for $T_L\equiv0$ (note that zero latent period gives an upper bound for $R_U$ for arbitrary latent period distribution), $\pi_R=1$, $\pi_T=0$, $p=1$ and $\xi=0.7$. The expression for $y^*$ in Eqn.~(\ref{eqn:ystar}) and $R_U$ in Eqn.~(\ref{eqn:RUexp}) feature infinite series, which, for calculation purposes, were truncated where further terms are zero to computational precision. 

We know that if $y$ is finite it is given by $y^*$. We see that $y^*$ increases monotonically before `blowing-up' for the first time (this occurs at around $\lambda=1.9876$). It then becomes negative, but, by definition, $y$ is non-negative and therefore cannot be given by $y^*$ and so must then be $\infty$ (and so also is $R_U$). Since $y$ is non-decreasing in $\lambda$, $y$ must be infinite always after this, even if $y^*$ is positive (which it would appear is possible). We assume then that $\lambda^*\approx1.9876$.

Figure \ref{fig:convergence} also shows a plot of $R_U$ varying with $\lambda$ (note that we are assuming $R_U$ is indeed finite for $\lambda\in[0,1.9876)$). We see that $R_U$ appears to asymptote at about $1.9876$, i.e.~as $y^*\to\infty$ for the first time. Similar behaviour was observed with other parameter values and latent period distributions.

\section{Numerical illustrations}
\label{sec:numerics}

Throughout this section we assume, without loss of generality, that the infectious period is of unit mean, i.e.~$\mbox{E}\left[T_I\right]=1$. All other time-related parameters (latent period and delay means, contact rates) should be thought of as being scaled accordingly.

\begin{figure}[h!]
\centering
\centerline{\includegraphics[width=1.15\textwidth]{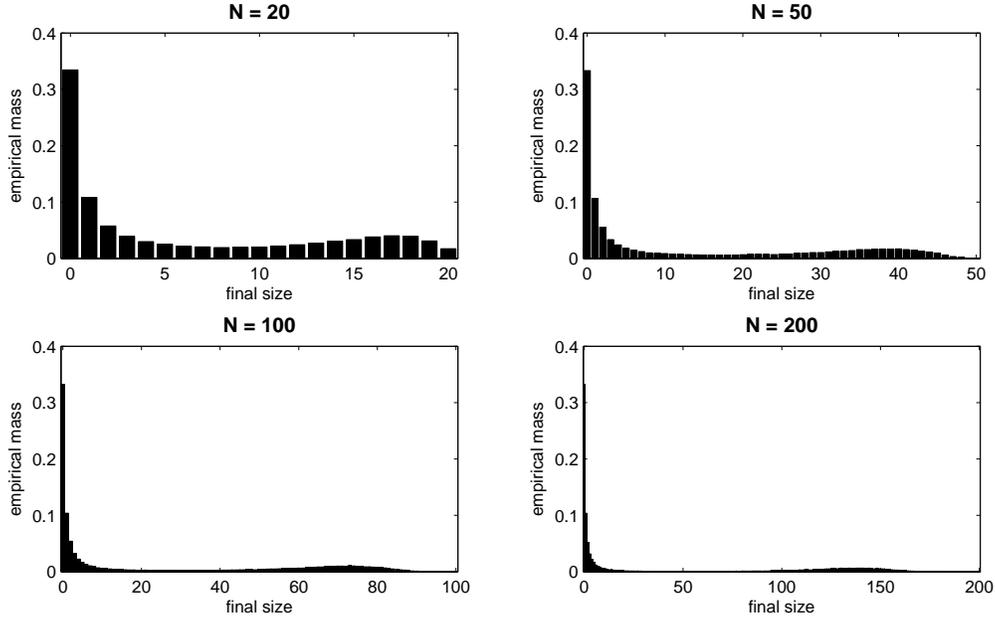}}
\caption{Final size distributions from 100,000 simulations, when $T_I\sim\mbox{Exp}(1)$, $T_D\sim\mbox{Exp}(1)$, $T_L\sim\mbox{Exp}(1)$, $\lambda=2$, $m=1$, $\pi_R=\pi_T=0.8$ and $p=0.5$.}
\label{fig:simulations-exp}
\end{figure}

\begin{figure}[h!]
\centering
\centerline{\includegraphics[width=1.15\textwidth]{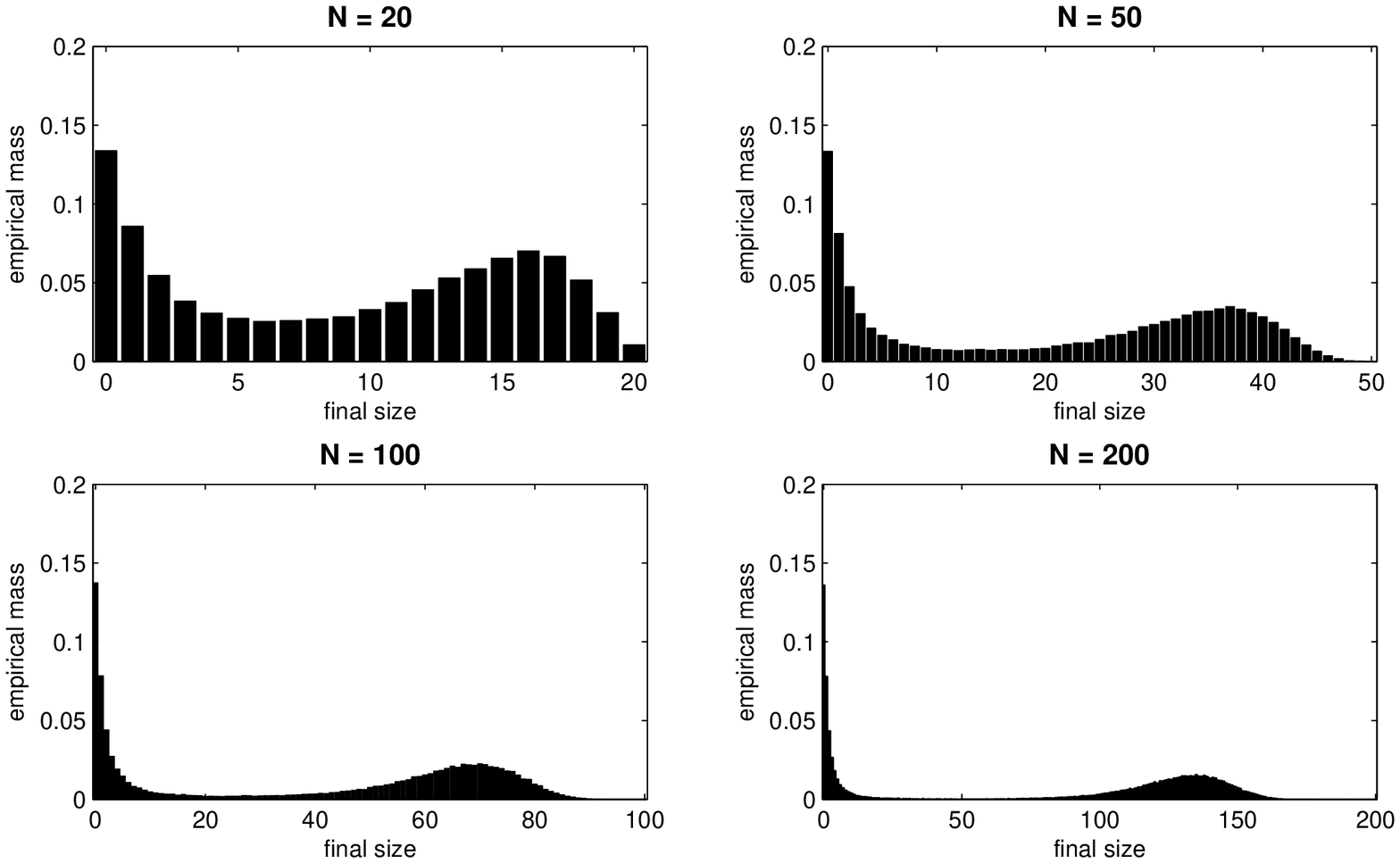}}
\caption{Final size distributions from 100,000 simulations, when $T_I\equiv1$, $T_D\sim\mbox{Exp}(1)$, $T_L\sim\mbox{Exp}(1)$, $\lambda=2$, $m=1$, $\pi_R=\pi_T=0.8$ and $p=0.5$.}
\label{fig:simulations-const}
\end{figure}

{\em The branching process approximation compares favourably with results from simulated epidemics.} While we consider the initially-susceptible population size, $N$, to be large to enable us to use branching process approximations to analyse our epidemic model, in real life $N$ is always finite, and often not `large' in a mathematical sense. Hence, it is of interest to examine how quickly the approximation becomes a valid description of the true epidemic model. Figures \ref{fig:simulations-exp} (Exponential infectious period) and \ref{fig:simulations-const} (Constant infectious period) show the final size (i.e.~total number of removals) distributions from $100,000$ simulations for $N=20$, $50$, $100$ and $200$ with $\mbox{E}\left[T_I\right]=1$, $\lambda=2$, $\pi_R=\pi_T=0.8$, $p=0.5$ and initial number of infectives $m=1$, when $T_L$ and $T_D$ both have $\mbox{Exp}(1)$ distributions. For the final size of an epidemic, we expect to see a bimodal distribution, with one mode corresponding to {\it minor outbreaks} (i.e.~a small proportion of the population is infected) and a second corresponding to {\it major outbreaks} (i.e.~a significant proportion is infected). We see this behaviour even for $N=20$, but there is not a clear distinction between minor and major outbreaks until about $N=200$.

\begin{figure}[h!]
\centering
\centerline{\includegraphics[width=1.1\textwidth]{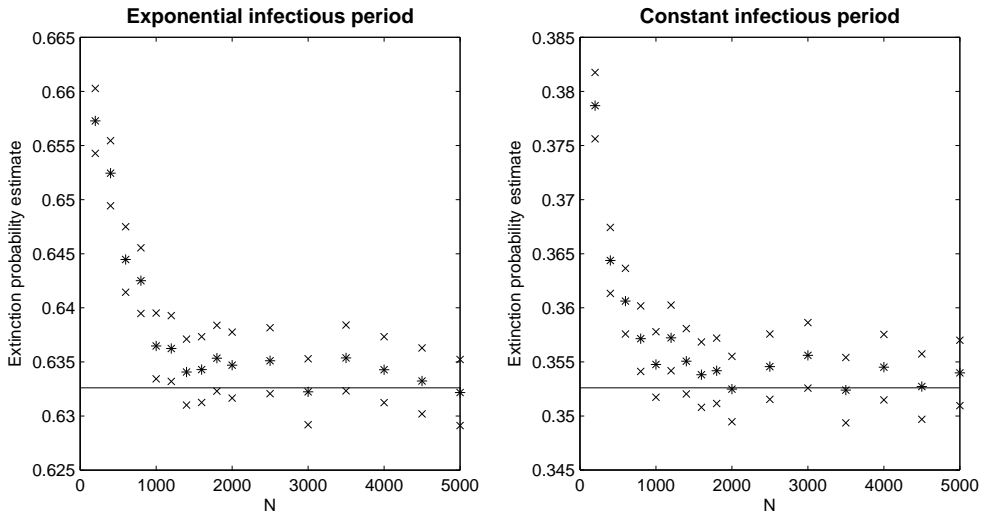}}
\caption{Estimated extinction probabilities (represented by asterisks) from $100,000$ simulations, when $\mbox{E}\left[T_I\right]=1$, $T_D\sim\mbox{Exp}(1)$, $T_L\sim\mbox{Exp}(1)$, $\lambda=2$, $m=1$, $\pi_R=\pi_T=0.8$ and $p=0.5$. Crosses represent two standard errors above and below the estimate and solid lines represent the true asymptotic extinction probabilities.}
\label{fig:estimatedpE}
\end{figure}

As $N\to\infty$, the proportion of minor outbreaks should tend to the theoretical extinction probability, $p_E$, for the branching process approximation. Thus, for $N$ large enough we estimate the extinction probability as the proportion of outbreaks that are minor. In Figure \ref{fig:estimatedpE} we plot these estimated extinction probabilities ($\hat{p}_E$), with confidence intervals given as $\hat{p}_E\pm2SE$ where the standard error is $SE=\left(\frac{\left(1-\hat{p}_E\right)\hat{p}_E}{n}\right)^{\frac{1}{2}}$ ($n=100,000$ is the number of simulations). These estimates are obtained by plotting the final size distributions of the simulations, determining a cut-off between major and minor outbreaks by sight, and defining our estimate as the proportion of outbreaks that are minor. Although we have not derived expressions for the extinction probability in these cases, our framework allows us to estimate it ($0.3526$ in the exponential case, $0.6326$ in the constant case) by simulating $R$ 100,000 times to obtain an empirical distribution for $R$; our estimate is then given as the solution of $\bar{H}(s)=s$, in $(0,1)$, where $\bar{H}(s)$ is the empirical probability generating function of $R$. In the both cases, the asymptotic extinction probability is consistently in the confidence interval for $N=1400$ and above. However, even for $N=200$ the estimates are reasonably close to the asymptotic value. In the limit $N\to\infty$, there are always infinitely many susceptibles, and thus infections do not affect a reduction in the number of susceptibles. For small $N$, infections noticeably reduce the number of susceptibles, and the probability of new infections, so there will be more minor outbreaks. This is why we see the asymptotic extinction probability being less than estimates obtained for small $N$.

\begin{figure}[h!]
\centering
\centerline{\includegraphics[width=1.15\textwidth]{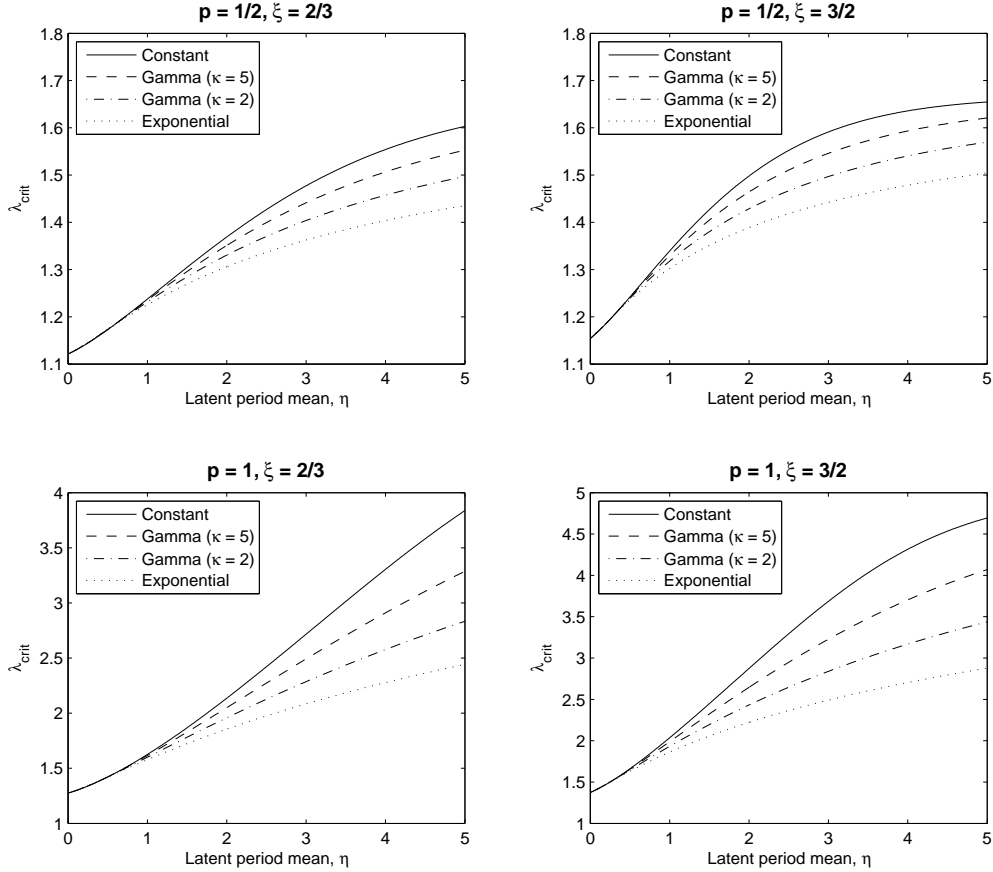}}
\caption{$\lambda_{crit}$ varying with mean latent period for different latent period distributions, when $T_I\sim\mbox{Exp}(1)$, $T_D\sim\mbox{Exp}(\xi)$ and $\pi_T=\pi_R=0.8$.}
\label{fig:lambdacrit-latent}
\end{figure}

{\em Latent period distribution can have a material effect on the spread of the epidemic, increasing with mean latent period.} Figure \ref{fig:lambdacrit-latent} shows the effect of latent period distribution choice in the model, plotting $\lambda_{crit}$, the critical contact rate (the contact rate which yields $R_U=1$) against the mean latent period for the latent period distributions: Exponential, Gamma (with shape parameter $\kappa=2,5$) and Constant. Note that if the latent period follows a Gamma distribution with mean $\mu$ and shape parameter $\kappa$ then $f_L(t)=(\kappa/\mu)^\kappa t^{\kappa-1}e^{-\kappa t/\mu}/\Gamma(\kappa)$ ($t>0$) where $\Gamma(\kappa)=\int_0^\infty t^{\kappa-1}e^{-t}\hspace{1mm}dt$ is the Gamma function, and $\phi_L(\theta)=(1+\mu\theta/\kappa)^\kappa$ ($\theta\geq0$). The infectious period has an $\mbox{Exp}(1)$ distribution and the delay has an $\mbox{Exp}(\xi)$ distribution (with $\xi$ chosen so that we have delays either typically shorter or typically longer than infectious periods). As the latent period mean increases, $\lambda_{crit}$ increases. We would expect this behaviour as a longer latent period increases the likelihood of an individual being traced, and further they experience less of their natural infectious period if they are traced. We can see that the effects of choosing a different latent period distribution increase as the mean latent period increases, as the mean delay decreases and as the naming probability increases. The difference between Exponential and Constant latent periods (Gamma distributions, with $\kappa=1$ and $\kappa=\infty$, respectively) is clearly distinct for given values of the other parameters.

\begin{figure}[h!]
\centerline{\includegraphics[width=1.15\textwidth]{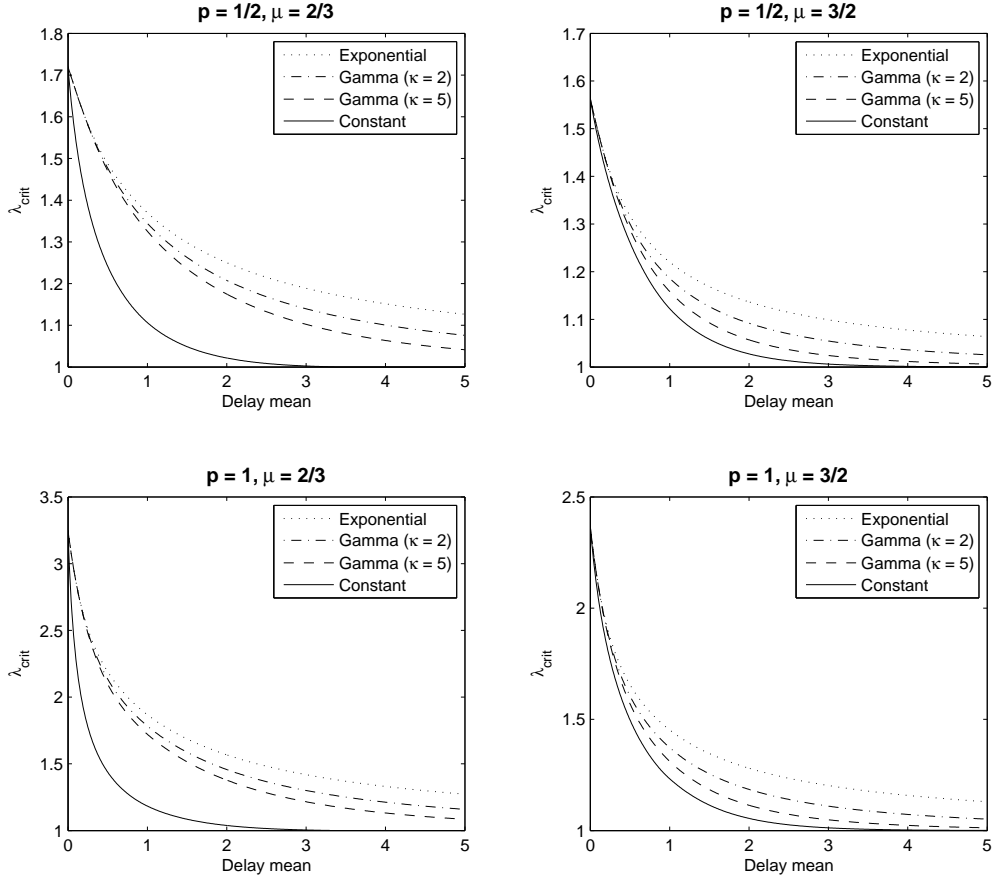}}
\caption{$\lambda_{crit}$ varying with delay mean for different delay distributions, when $T_I\equiv1$, $T_L\sim\mbox{Exp}(\mu)$, $\pi_R=1$ and $\pi_T=0$.}
\label{fig:lambdacrit-delay}
\end{figure}

{\em Delay distribution can have a material effect on the spread of the epidemic.} Figure \ref{fig:lambdacrit-delay} shows the effect of delay distribution choice in the model, plotting $\lambda_{crit}$ against the mean delay for the delay distributions: Exponential, Gamma (with shape parameter $\kappa=2,5$) and Constant. The infectious period is constant ($T_I\equiv1$) and the latent period has an $\mbox{Exp}(\mu)$ distribution (with $\mu$ chosen so that we have latent periods either typically shorter or typically longer than infectious periods). As delay mean increases, $\lambda_{crit}$ decreases. We would expect this as a longer delay means an individual is less likely to be traced, while if they are they serve more of their natural infectious period. We can see that the effects of choosing a different delay distribution increase as the delay mean increases from zero (though as the delay mean tends to infinity, $\lambda_{crit}\to1$, irrespective of the exact distribution), and as the naming probability increases. It would appear that as the shape parameter $\kappa$ of the Gamma distribution increases, there is much slower convergence of $\lambda_{crit}$ to that in the constant case for smaller $\mu$ (i.e.~longer latent periods). Overall an exponentially-distributed delay (and so higher variance) is better for control, which can be explained by (a) having a median lower than the mean, so typically delays are shorter, and (b) the effect of allowing for longer delays being mitigated by the fact that once a delay is long enough that an individual is not traced, it does not matter how long that delay actually is.

{\em Tracing delays of the same length and of independent length for individuals with the same infector yield similar results.} While we have assumed that sibling units experience independent delays, it has been seen that $R_U$ is unchanged by assuming instead that sibling units experience the same delay. However, the probability of extinction is not unchanged, so it is of interest to see how much it differs between the two delay assumptions.

Table \ref{table:mutvsinddelays} shows how the extinction probabilities compare for independent and mutual delays when $\pi_R=\pi_T=0.8$, $p=1$ and $T_I\sim\mbox{Exp}(1)$. The extinction probabilities are estimated as previously, although for longer delays, $R$ is infinite with non-zero probability; to circumvent this we assume that if the number of offspring is sufficiently large (here we consider this to mean at least $100$), then it is infinite. The results are shown for $\lambda=1.5$ and $T_L\equiv0$, $\lambda=2.5$ and $T_L\equiv0$, and $\lambda=1.5$ and $T_L\sim\mbox{Exp}(1)$.

We see that there is no dramatic difference between the two delay assumptions, even as we increase $\lambda$ or introduce a latent period. There is a slight tendency for higher extinction probabilities with mutual delays, which concurs with the findings of Ball {\em et al.~}\cite{ball-knock-oneill}, that dependencies between sibling units increase the extinction probability. Similar results were obtained when different assumptions were tested, such as constant infectious periods, constant latent periods or different values for $\pi_R$, $\pi_T$ and $p$.

\begin{table}
\begin{tabular}{|c|c|c|c|c|c|c|}
\hline
Delay & \multicolumn{2}{|c|}{$\lambda=1.5$, $T_L\sim\mbox{Exp}(1)$} & \multicolumn{2}{|c|}{$\lambda=1.5$, $T_L\equiv0$} & \multicolumn{2}{|c|}{$\lambda=2.5$, $T_L\equiv0$}\\
\cline{2-7}
mean & Independent & Mutual & Independent & Mutual & Independent & Mutual\\
\hline
$0.0$ & $1.0000$ & $1.0000$ & $1.0000$ & $1.0000$ & $0.8300$ & $0.8270$\\
$0.5$ & $1.0000$ & $1.0000$ & $0.9821$ & $0.9861$ & $0.5166$ & $0.5296$\\
$1.0$ & $1.0000$ & $1.0000$ & $0.8438$ & $0.8463$ & $0.4632$ & $0.4709$\\
$1.5$ & $0.9999$ & $0.9999$ & $0.7865$ & $0.7884$ & $0.4434$ & $0.4493$\\
$2.0$ & $0.9718$ & $0.9733$ & $0.7533$ & $0.7586$ & $0.4314$ & $0.4366$\\
$2.5$ & $0.9090$ & $0.9191$ & $0.7366$ & $0.7426$ & $0.4239$ & $0.4291$\\
$3.0$ & $0.8723$ & $0.8781$ & $0.7280$ & $0.7286$ & $0.4206$ & $0.4223$\\
$3.5$ & $0.8431$ & $0.8504$ & $0.7184$ & $0.7175$ & $0.4171$ & $0.4229$\\
$4.0$ & $0.8232$ & $0.8288$ & $0.7111$ & $0.7119$ & $0.4159$ & $0.4199$\\
$4.5$ & $0.8041$ & $0.8121$ & $0.7039$ & $0.7089$ & $0.4143$ & $0.4162$\\
$5.0$ & $0.7875$ & $0.7944$ & $0.7032$ & $0.7021$ & $0.4130$ & $0.4149$\\
\hline
\end{tabular}
\caption{Extinction probability estimated from simulations for varying delay mean ($\frac{1}{\xi}$) for both mutual and independent delays when $T_I\sim\mbox{Exp}(1)$, $T_D\sim\mbox{Exp}(\xi)$, $\pi_R=\pi_T=0.8$ and $p=1$.}
\label{table:mutvsinddelays}
\end{table}

\begin{figure}[h!]
\centerline{\includegraphics[width=1.15\textwidth]{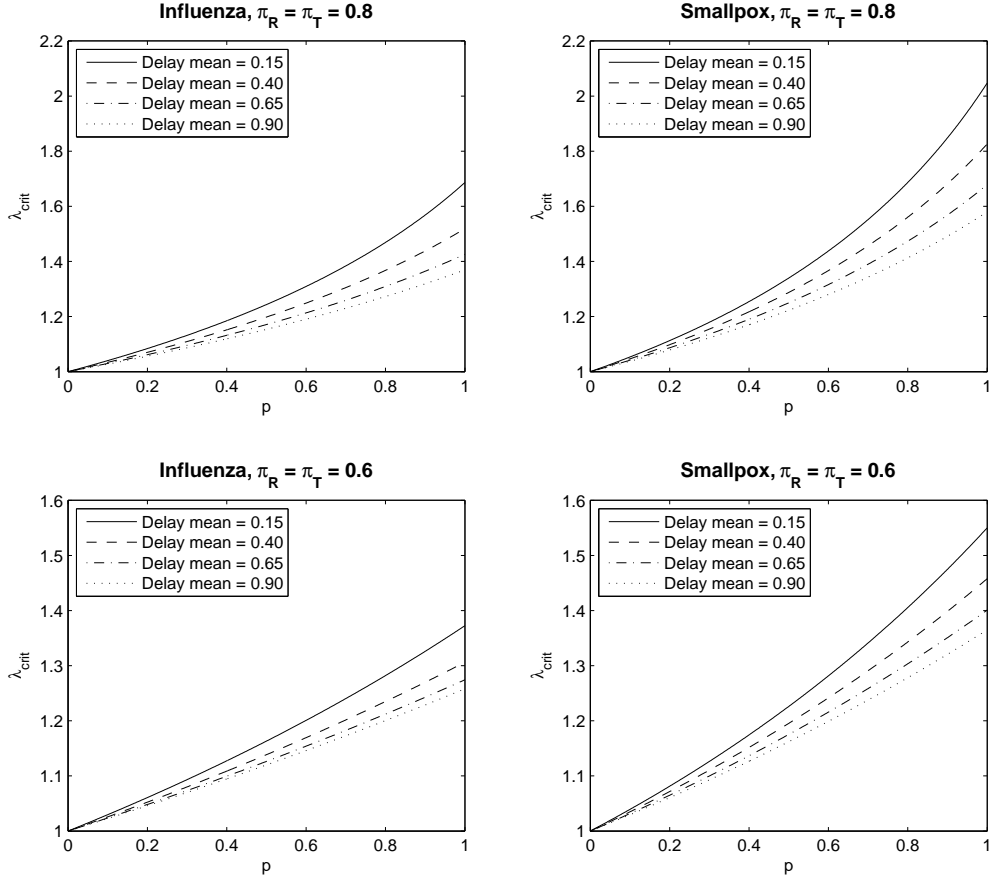}}
\caption{$\lambda_{crit}$ varying with $p$ for different delay means ($1/\xi$), when $T_I\sim\mbox{Exp}(1)$, $T_D\sim\mbox{Exp}(\xi)$, $T_L\sim\mbox{Exp}(1/0.10)$ (influenza) and $T_L\sim\mbox{Exp}(1/0.58)$ (smallpox).}
\label{fig:influenzasmallpox}
\end{figure}

{\em Contact tracing better controls diseases with longer latent periods and more infected individuals able to be interviewed.} We consider now an application to some infectious diseases. Table 2 of Klinkenberg {\em et al.}~\cite{klink-fraser-heest} suggests the mean latent periods of influenza and smallpox are about $0.10$ and $0.58$ (relative to unit mean infectious period), respectively. In Figure \ref{fig:influenzasmallpox}, we see how $\lambda_{crit}$ varies with $p$ and delay mean for these two latent period means. The contact tracing is more effective for smallpox, with its longer latent period. The effect of reducing delays becomes greater as the naming probability increases, and vice versa.

\section{Discussion}
\label{sec:discussion}

In this paper we restricted attention to forward tracing, that is that tracing only occurs by infectors naming their infectees. This allows for simpler analysis since each individual depends directly upon only one other individual (their infector). We have not considered backward tracing (infectees naming their infectors) which is more difficult to analyse. Hethcote \& Yorke \cite{hethcote-yorke} suggest tracing from infectee to infector is more effective with regards to heterogeneous mixing for gonorrhea, while M\"{u}ller {\em et al.}~\cite{muller-kretzschmar-dietz} suggest the opposite for homogeneous mixing (as we have here). With latent periods and delays, we would expect the effect of backward tracing to be reduced since individuals are more likely to be removed before they can be traced via an infectee.

We have considered a population that mixes homogeneously, however incorporating population structure can enhance contract tracing. Eames \& Keeling \cite{eames-keeling} model contact tracing on a network, where triangles allow infectees of the same infector to be traced via one another. Shaban {\em et al.}~\cite{shabanetal} also model contact tracing on a network (without clustering), with a vaccination model that necessarily means the contact tracing targets susceptible neighbours at high-risk of infection. It would be beneficial to consider contact tracing that targets intervention at both a  named infected individual and the individuals they are more likely to infect. We will consider this via a contact tracing model for an epidemic spreading amongst a community of households in a separate paper.

For the purposes of controlling an epidemic, it may be easier in practice to reduce the length of delays than increase the naming probability. Lemma 1 gives an indication that short delays are important. However, numerical results show that, for the same delay mean, higher variance of delays is better for control as more individuals are traced while latent or early on in their infectious period than for a low variance.

\section*{Acknowledgements}
Edward Knock was supported by an EPSRC Doctoral Training grant.

\appendix
\section{Proof of Lemma 1}
\label{app:RUinfty}
\begin{proof}

Suppose first that $\mbox{P}\left(T_D>T_L\right)=0$. The case where $T_D\equiv T_L$ essentially reduces to the case with no latent periods or delays, for which, when $\pi_R=\pi_T=1$ (i.e.~everyone is interviewed), $R_U=\frac{1-p}{p}\left(e^{\lambda p/\gamma}-1\right)$ (Eqn.~(5.1) in Ball {\em et al.}~\cite{ball-knock-oneill}). For more general $\left(\pi_R,\pi_T\right)$ in this case, fewer individuals are interviewed, so the number of immediate type-$k$ ($k=1,2,\ldots$) descendants of a typical unnamed individual is stochastically smaller, but this unnamed individual and some of these descendants may not be interviewed and so have unnamed offspring at rate $\lambda$ instead of $\lambda(1-p)$. Hence, $R_U<\frac{1}{1-p}\times\frac{1-p}{p}\left(e^{\lambda p/\gamma}-1\right)$, so $R_U<\infty$ if $\lambda<\infty$. If a named individual's latent period cannot be shorter than the delay in tracing them, then active lifetimes will be at not longer than in the $T_D\equiv T_L$ case, whence $R$ is stochastically smaller. Hence, if $\mbox{P}\left(T_D>T_L\right)=0$, then $R_U<\infty$.

Suppose now that $\mbox{P}\left(T_D>T_L\right)>0$. First we consider a lower bound for $R_U$ in order to show that it is infinite for $\lambda$ sufficiently large. Recall that in the MBDP, $R$ is the total number of unnamed immediate offspring of (a) an unnamed individual and (b) of all named descendants of this unnamed individual, for whom the unnamed individual is their nearest unnamed ancestor. Consider $R^-$, a lower bound for $R$ in which instead of counting the unnamed offspring of {\it all} named descendants as described, we count only (a) the immediate unnamed offspring of the unnamed individual and (b) unnamed individuals who are separated in the family tree from the unnamed individual by only named ancestors that are asked to name their offspring and have a natural lifetime in the interval $(0,\epsilon)$ ($\epsilon>0$) and an associated delay at least $\epsilon$ time units greater than their latent period (note that, since $\mbox{P}\left(T_D>T_L\right)>0$ there does exist some $\epsilon>0$ such that $\mbox{P}\left(T_D>T_L+\epsilon\right)>0$). Clearly then, $R^-\leq R$, and if $R^-_U=\ez\left[R^-\right]$, then $R^-_U\leq R_U$.

Named individuals who have a natural lifetime in the interval $(0,\epsilon)$ and who experience a delay at least $\epsilon$ time units greater than their latent period, must end their lifetime before the delay ends, and hence are necessarily untraced. Therefore they are asked to name their offspring with probability $\pi_R$, in which case they produce named offspring at rate $\lambda p$ over the length of their lifetime, and the probability that a typical named offspring has a natural lifetime in the interval $(0,\epsilon)$ and  experience a delay at least $\epsilon$ time units greater than their latent period is $\left(1-e^{-\gamma\epsilon}\right)\pz\left(T_D>T_L+\epsilon\right)$, and so these particular named offspring are produced at rate $\lambda p\left(1-e^{-\gamma\epsilon}\right)\pz\left(T_D>T_L+\epsilon\right)$. Such named individuals behave independently of one another, and hence are described by a branching process. Now
\begin{eqnarray*}
\ez\left[T_I\left|0<T_I<\epsilon\right.\right]=\left[1-e^{-\gamma\epsilon}(\gamma\epsilon+1)\right]/\left[\gamma\left(1-e^{-\gamma\epsilon}\right)\right],
\end{eqnarray*}
so the expected number of offspring of a typical individual in this branching process is
\begin{eqnarray*}
\frac{\lambda p\pi_R}{\gamma}\left[1-e^{-\gamma\epsilon}(\gamma\epsilon+1)\right]\pz\left(T_D>T_L+1\right),
\end{eqnarray*}
whence the total progeny is infinite with positive probability if
\begin{eqnarray*}
\lambda > \left\{\frac{p\pi_R}{\gamma}\left[1-e^{-\gamma\epsilon}(\gamma\epsilon+1)\right]\pz\left(T_D>T_L+1\right)\right\}^{-1}.
\end{eqnarray*}

Now since these individuals have positive lifetimes and in the MBDP give birth to unnamed individuals at rate $\lambda(1-p)$, for large enough $\lambda$, $R^-_U=\infty$ and so $R_U=\infty$. Now, $R_U$ is non-decreasing in $\lambda$, so if we define $\lambda^*=\inf\left\{\lambda:R_U=\infty\right\}$, then for all $\lambda>\lambda^*$, $R_U$ is infinite. However, if $\lambda\in\left[0,\lambda^*\right)$ then $R_U$ is finite. Since the total progeny of the birth-death process without tracing has finite mean if $\lambda/\gamma<1$, then it must be the case that $R_U<\infty$ if $\lambda/\gamma<1$, so $\lambda^*\geq\gamma$.
\end{proof}

\section{Details of the derivation of $L_k(\theta)$ in Section \ref{sec:exp}}
\label{app:Lkthetadetails}

It follows from Eqn.~(\ref{eqn:hk1}) that
\begin{align*}
L_k(\theta) & = \int_0^\infty e^{-\theta t}\left(\lambda\left(1-\pi_Rp\right) t+\lambda\pi_Rptg_k(t)\right)\hspace{1mm}dt\\
            & = \frac{\lambda\left(1-\pi_Rp\right)}{\theta^2} + \lambda\pi_Rp\int_0^\infty e^{-\theta t}tg_k(t)\hspace{1mm}dt.
\end{align*}

For $\theta\geq0$, let $\psi(\theta)=\int_0^\infty e^{-\theta w}\chi(w)\hspace{1mm} dw$, then
\begin{align*}
\psi(\theta)& =\int_0^\infty e^{-\theta v}\{\lambda(1-\pi_T) ve^{-\gamma v}+\pi_T e^{-\gamma v}h_{k-1}(v)\\
&\hspace{2.5cm}+\int_{u=0}^v\gamma e^{-\gamma u}\left(\lambda\left(1-\pi_R\right)u+\pi_Rh_{k-1}(u)\right)\hspace{1mm}du\}dv\\
& = \frac{\lambda(1-\pi_T)}{(\theta+\gamma)^2}+\frac{\lambda\gamma(1-\pi_R)}{\theta(\theta+\gamma)^2}+\frac{\pi_T\theta+\pi_R\gamma}{\theta}L_{k-1}(\theta+\gamma).
\end{align*}

Noting that
\begin{align*}
&   \int_{t=0}^\infty e^{-\theta t}\int_{t_D=0}^\infty\xi e^{-\xi t_D}\int_{t_L=0}^{t_D} \int_{v=t_D-t_L}^{t+t_D-t_L}\chi(v)f_L\left(t_L\right) \hspace{1mm} dv\hspace{1mm} dt_L\hspace{1mm} dt_D \hspace{1mm} dt\displaybreak[0]\\
& = \int_{t=0}^\infty e^{-\theta t}\int_{t_L=0}^\infty\int_{t_D=t_L}^{\infty} \xi e^{-\xi t_D}\int_{v=t_D-t_L}^{t+t_D-t_L}\chi(v)f_L\left(t_L\right)\hspace{1mm} dv \hspace{1mm}dt_D\hspace{1mm} dt_L\hspace{1mm} dt \displaybreak[0]\\
& = \int_{t=0}^\infty e^{-\theta t}\int_{t_L=0}^\infty e^{-\xi t_L}f_L\left(t_L\right)\int_{t'_D=0}^{\infty} \xi e^{-\xi t'_D}\int_{v=t'_D}^{t+t'_D}\chi(v) \hspace{1mm} dv \hspace{1mm}dt'_D\hspace{1mm} dt_L\hspace{1mm} dt \displaybreak[0]\\
& = \phi_L(\xi)\int_{t=0}^\infty e^{-\theta t}\int_{t'_D=0}^{\infty} \xi e^{-\xi t'_D}\int_{v=t'_D}^{t+t'_D}\chi(v) \hspace{1mm} dv \hspace{1mm}dt'_D \hspace{1mm}dt \displaybreak[0]\\
& = \phi_L(\xi)\int_{t'_D=0}^{\infty} \xi e^{-\xi t'_D}\int_{v=t'_D}^{\infty}\chi(v)\int_{t=v-t'_D}^\infty e^{-\theta t}\hspace{1mm} dt\hspace{1mm} dv \hspace{1mm}dt'_D \displaybreak[0]\\
& = \frac{\phi_L(\xi)}{\theta}\int_{v=0}^{\infty}e^{-\theta v}\chi(v)\int_{t'_D=0}^{v} \xi e^{-(\xi-\theta) t'_D}\hspace{1mm} dt'_D\hspace{1mm} dv \displaybreak[0]\\
& = \left\{\begin{array}{ll}
\frac{\xi\phi_L(\xi)}{\theta(\xi-\theta)}\int_{v=0}^{\infty}\left(e^{-\theta v}-e^{-\xi v}\right)\chi(v) \hspace{1mm}dv & \mbox{if $\theta\neq\xi$}\\
\phi_L(\xi)\int_{v=0}^{\infty}ve^{-\xi v}\chi(v)\hspace{1mm} dv & \mbox{if $\theta=\xi$}\end{array}\right. \displaybreak[0]\\
& = \left\{\begin{array}{ll}
\frac{\xi\phi_L(\xi)}{\theta(\xi-\theta)}\left(\psi(\theta)-\psi(\xi)\right) & \mbox{if $\theta\neq\xi$}\\
\phi_L(\xi)\psi'(\xi) & \mbox{if $\theta=\xi$},\end{array}\right.
\end{align*}
and further that
\begin{align*}
&   \int_{t=0}^\infty e^{-\theta t}\int_{t_D=0}^\infty\xi e^{-\xi t_D}\int_{t_L=t_D}^{t_D+t} \int_{v=0}^{t+t_D-t_L}\chi(v)f_L\left(t_L\right)\hspace{1mm}  dv\hspace{1mm} dt_L\hspace{1mm} dt_D\hspace{1mm} dt\displaybreak[0]\\
& = \int_{t_D=0}^\infty\xi e^{-\xi t_D}\int_{t_L=t_D}^\infty\int_{t'=0}^\infty e^{-\theta (t'-t_D+t_L)}\int_{v=0}^{t'}\chi(v)f_L\left(t_L\right) \hspace{1mm} dv \hspace{1mm}dt'\hspace{1mm} dt_L\hspace{1mm} dt_D \displaybreak[0]\\
& = \int_{t_D=0}^\infty\xi e^{-\xi t_D}\int_{t_L=t_D}^\infty f_L\left(t_L\right)\int_{v=0}^{\infty}\chi(v)\int_{t'=v}^\infty e^{-\theta (t'-t_D+t_L)} \hspace{1mm} dt' \hspace{1mm} dv \hspace{1mm} dt_L\hspace{1mm} dt_D \displaybreak[0]\\
& = \frac{1}{\theta} \int_{t_L=0}^\infty e^{-\theta t_L}f_L\left(t_L\right) \int_{t_D=0}^{t_L}\xi e^{-(\xi-\theta)t_D} \int_{v=0}^{\infty}e^{-\theta v}\chi(v) \hspace{1mm} dv \hspace{1mm} dt_D \hspace{1mm} dt_L \displaybreak[0]\\
& = \frac{1}{\theta} \int_{t_L=0}^\infty e^{-\theta t_L} f_L\left(t_L\right)\int_{t_D=0}^{t_L}\xi e^{-(\xi-\theta)t_D} \psi(\theta) \hspace{1mm} dt_D \hspace{1mm}dt_L \displaybreak[0]\\
& = \left\{\begin{array}{ll}
\frac{\xi}{\theta(\xi-\theta)} \int_{t_L=0}^\infty\left(e^{-\theta t_L}-e^{-\xi t_L}\right)f_L\left(t_L\right)\psi(\theta) \hspace{1mm} dt_L & \mbox{if $\theta\neq\xi$} \\
\int_{t_L=0}^\infty t_Le^{-\xi t_L} f_L\left(t_L\right)\psi(\xi) \hspace{1mm}dt_L & \mbox{if $\theta=\xi$} \end{array}\right.\displaybreak[0]\\
& = \left\{\begin{array}{ll}
\frac{\xi\left(\phi_L(\theta)-\phi_L(\xi)\right)}{\theta(\xi-\theta)} \psi(\theta) & \mbox{if $\theta\neq\xi$} \\
-\phi'_L(\xi) \psi(\xi) & \mbox{if $\theta=\xi$},\end{array}\right.\displaybreak[0]
\end{align*}
it follows that, for $\theta\neq\xi$,
\begin{align*}
L_k(\theta) & = \frac{\lambda(1-\pi_Rp)}{\theta^2}+ \frac{\lambda \pi_Rp\xi\left(\phi_L(\theta)-\phi_L(\xi)\right)}{\theta(\xi-\theta)} \psi(\theta)\displaybreak[0]\\
            & \quad + \frac{\lambda\pi_R p\xi\phi_L(\xi)}{\theta(\xi-\theta)}\left(\psi(\theta)-\psi(\xi)\right)\displaybreak[0]\\
            & = \frac{\lambda(1-\pi_Rp)}{\theta^2}+ \frac{\lambda\pi_R p\xi}{\theta(\xi-\theta)} \left\{\phi_L(\theta)\psi(\theta)-\phi_L(\xi)\psi(\xi)\right\},
\end{align*}
while
\begin{align*}
L_k(\xi) & = \frac{\lambda(1-\pi_Rp)}{\xi^2}+\lambda\pi_Rp\left\{\phi_L(\xi)\psi'(\xi)-\phi'_L(\xi)\psi(\xi)\right\}\displaybreak[0]\\
         & = \frac{\lambda(1-\pi_Rp)}{\xi^2}+\lambda\pi_Rp\lim_{\theta\to\xi}\left\{\phi_L(\xi)\frac{\psi(\theta)-\psi(\xi)}{\theta-\xi}-\frac{\phi_L(\theta)-\phi_L(\xi)}{\theta-\xi}\psi(\xi)\right\}\displaybreak[0]\\
         & = \lim_{\theta\to\xi}L_k(\theta),
\end{align*}
thus we can obtain $L_k(\theta)$ when $\theta=\xi$ by taking the limit $\theta\to\xi$. From now on we assume that $\theta\neq\xi$, and that $\xi$ is non-integer. If $\xi$ is an integer, we can obtain $R_U$ by considering the limit as $\xi$ approaches this integer value.

Setting $\gamma=1$, without loss of generality, yields
\begin{align*}
\psi(\theta)& = \frac{\lambda(1-\pi_T)}{(\theta+1)^2}+\frac{\lambda(1-\pi_R)}{\theta(\theta+1)^2}+\frac{\pi_T\theta+\pi_R}{\theta}L_{k-1}(\theta+1).
\end{align*}

\section{Details of the derivation of $x_k(\theta)$ in Section \ref{sec:exp}}
\label{app:xkdetails}

If $x_{k,j}(\theta)=L_k(j+\theta)$ ($\theta\geq0$, $k=0,1,\ldots$, $j=1,2,\ldots$), then $x_{0,j}(\theta)=\frac{\lambda(1-\pi_Rp)}{(j+\theta)^2}$ and, for $k=1,2,\ldots$ and $j=1,2,\ldots$,
\begin{eqnarray}
x_{k,j}(\theta)=a_j(\theta)-\rho_j(\theta)x_{k-1,1}(\xi)+b_j(\theta)x_{k-1,j+1}(\theta).
\label{xkjeqn}
\end{eqnarray}

Suppose that, for $1\leq\kappa\leq k-1$,
\begin{align}
\begin{split}
x_{\kappa,k-\kappa+1}(\theta) & = \frac{c_k(\theta)}{c_{k-\kappa}(\theta)}x_{0,k+1}(\theta)\\
							  & \quad +\sum_{i=k-\kappa}^{k-1}\frac{c_i(\theta)}{c_{k-\kappa}(\theta)} \left(a_{i+1}(\theta)-\rho_{i+1}(\theta)x_{k-i-1,1}(\xi)\right),\\
\end{split}
\label{eqn:xkappa}
\end{align}
then using Eqn.~(\ref{xkjeqn}),
\begin{align*}
x_{\kappa+1,k-\kappa}(\theta) & = a_{k-\kappa}(\theta)-\rho_{k-\kappa}(\theta)x_{\kappa,1}(\xi)+b_{k-\kappa}(\theta)x_{\kappa,k-\kappa+1}(\theta) \displaybreak[0]\\
                      		  & = a_{k-\kappa}(\theta)-\rho_{k-\kappa}(\theta)x_{\kappa,1}(\xi)+b_{k-\kappa}(\theta)\left(\frac{c_k(\theta)}{c_{k-\kappa}(\theta)}x_{0,k+1}(\theta)\right.\\
                              & \quad \left.\hspace{2cm}+\sum_{i=k-\kappa}^{k-1}\frac{c_i(\theta)}{c_{k-\kappa}(\theta)} \left(a_{i+1}(\theta)-\rho_{i+1}(\theta)x_{k-i-1,1}(\xi)\right)\right) \displaybreak[0]\\
                              & = a_{k-\kappa}(\theta)-\rho_{k-\kappa}(\theta)x_{\kappa,1}(\xi)+\frac{c_k(\theta)}{c_{k-\kappa-1}(\theta)}x_{0,k+1}(\theta)\\
                              & \quad \hspace{1cm}+\sum_{i=k-\kappa}^{k-1}\frac{c_i(\theta)}{c_{k-\kappa-1}(\theta)} \left(a_{i+1}(\theta)-\rho_{i+1}(\theta)x_{k-i-1,1}(\xi)\right) \displaybreak[0]\\
                              & = \frac{c_k(\theta)}{c_{k-\kappa-1}(\theta)}x_{0,k+1}(\theta)\\
                              & \quad +\sum_{i=k-\kappa-1}^{k-1}\frac{c_i(\theta)}{c_{k-\kappa-1}(\theta)} \left(a_{i+1}(\theta)-\rho_{i+1}(\theta)x_{k-i-1,1}(\xi)\right),
\end{align*}
hence Eqn.~(\ref{eqn:xkappa}) is true for $1\leq\kappa\leq k$, and setting $\kappa=k$ yields
\begin{eqnarray}
x_{k,1}(\theta) = \sum_{i=0}^{k-1}c_i(\theta) \left(a_{i+1}(\theta)-\rho_{i+1}(\theta)x_{k-i-1,1}(\xi)\right)+c_k(\theta)x_{0,k+1}(\theta).
\label{eqn:xk1}
\end{eqnarray}

\section{Details of the derivation of $x(\theta)$ in Section \ref{sec:exp}}
\label{app:xthetadetails}

Note that $y=\lim_{k\to\infty}x_{k,1}(\xi)$ and $R_U=\lim_{k\to\infty}x_{k,1}(0)$. Since $x_{k,1}(\xi)\leq x_{k,1}(0)$ for all $k$, if $y=\infty$ then $R_U=\infty$.

Suppose then that $y<\infty$. Since $\left|\rho_{i}(\theta)\right|\leq\frac{\lambda\pi_R p(\pi_T\xi+\pi_R)}{|\xi-[\xi]|}$ and $\left|c_i(\theta)\right|\leq\frac{1}{i!}\left(\frac{\lambda\pi_R p\xi(\pi_T+\pi_R)}{|\xi-[\xi]|}\right)^i$, we have that
\begin{align*}
\sum_{i=0}^{\infty}\left|c_i(\theta)\rho_{i+1}(\theta)\right| & \leq \frac{\lambda\pi_R p(\pi_T\xi+\pi_R)}{|\xi-[\xi]|}\sum_{i=0}^{\infty}\frac{1}{i!}\left(\frac{\lambda\pi_R p\xi(\pi_T+\pi_R)}{|\xi-[\xi]|}\right)^i\\
                                                         & = \frac{\lambda\pi_R p(\pi_T\xi+\pi_R)}{|\xi-[\xi]|}e^{\frac{\lambda\pi_R p\xi(\pi_T+\pi_R)}{|\xi-[\xi]|}},
\end{align*}
and hence $\sum_{i=0}^{\infty}c_i(\theta)\rho_{i+1}(\theta)$ is convergent, and by the dominated convergence theorem for series, as $k\to\infty$,
\begin{align*}
\sum_{i=0}^{\infty}c_i(\theta)\rho_{i+1}(\theta)x_{k-i-1,1}(\xi)\to y\sum_{i=0}^{\infty}c_i(\theta)\rho_{i+1}(\theta).
\end{align*}
Further, $c_k(\theta)x_{0,k+1}(\theta)\to0$ as $k\to\infty$ and, since $\left|a_i(\theta)\right|\leq\lambda(1-\pi_Rp)+\frac{2\lambda^2\pi_Rp\xi}{|\xi-[\xi]|}$,
\begin{align*}
\sum_{i=0}^{\infty}\left|c_i(\theta)a_{i+1}(\theta)\right| & \leq \left(\lambda(1-\pi_Rp)+\frac{2\lambda^2\pi_Rp\xi}{|\xi-[\xi]|}\right) \sum_{i=0}^{\infty} \frac{1}{i!}\left(\frac{\lambda\pi_R p\xi(\pi_T+\pi_R)}{|\xi-[\xi]|}\right)^i\\
                                           & = \left(\lambda(1-\pi_Rp)+\frac{2\lambda^2\pi_Rp\xi}{|\xi-[\xi]|}\right) e^{\frac{\lambda\pi_R p\xi(\pi_T+\pi_R)}{|\xi-[\xi]|}},
\end{align*}
so $\sum_{i=0}^{\infty}c_i(\theta)a_{i+1}(\theta)$ is convergent.

Hence, letting $k\to\infty$ in Eqn.~(\ref{eqn:xk1}),
\begin{eqnarray*}
\lim_{k\to\infty}x_{k,1}(\theta) = \sum_{i=0}^{\infty}c_i(\theta) \left(a_{i+1}(\theta)-\rho_{i+1}(\theta)y\right),
\end{eqnarray*}
which is convergent, so $R_U<\infty$. Therefore, $y<\infty$ if and only if $R_U<\infty$.


\bibliographystyle{elsarticle-num}



\end{document}